# A growth model in multiple dimensions and the height of a random partial order

Timo Seppäläinen[1],*

*University of Wisconsin-Madison*

**Abstract:** We introduce a model of a randomly growing interface in multi-dimensional Euclidean space. The growth model incorporates a random order model as an ingredient of its graphical construction, in a way that replicates the connection between the planar increasing sequences model and the one-dimensional Hammersley process. We prove a hydrodynamic limit for the height process, and a limit which says that certain perturbations of the random surface follow the characteristics of the macroscopic equation. By virtue of the space-time Poissonian construction, we know the macroscopic velocity function explicitly up to a constant factor.

## 1. Introduction

We introduce a model of a randomly growing interface, whose construction involves the height of a random partial order. The interface is defined by a height function on $d$-dimensional Euclidean space, and the related model of random order is in $d + 1$ dimensional space-time. Our goal is to emulate in higher dimensions the fruitful relationship between the one-dimensional Hammersley process and the model of increasing sequences among planar Poisson points. The connection between Hammersley's process and increasing sequences was suggested in Hammersley's paper [15], first utilized by Aldous and Diaconis [1], and subsequently in papers [21, 22, 25, 27]. A review of the use of Hammersley's process to study increasing sequences appeared in [14], and of the wider mathematical context in [2]. The study of higher dimensional random orders was started by Winkler [30].

The interface process we introduce is defined through a graphical representation which utilizes a homogeneous space-time Poisson point process, and in particular the *heights of the partial orders* among the Poisson points in space-time rectangles. This definition suggests a natural infinitesimal description, which we verify in a sense. After defining the process, we prove a hydrodynamic limit for the height function. This proceeds in a familiar way, by the path level variational formulation. The deterministic limiting height is the solution of a Hamilton-Jacobi equation given by a Hopf-Lax formula.

Next we use this process to prove a limit that in a way generalizes the law of large numbers of a second class particle in one-dimensional systems. In interacting particle systems, a second class particle is the location $X(t)$ of the unique discrepancy between two coupled systems that initially differ by exactly one particle (see

---

*Research partially supported by NSF Grants DMS-01-26775 and DMS-04-02231.
[1]Mathematics Department, University of Wisconsin-Madison, Madison, WI 53706-1388, USA
e-mail: seppalai@math.wisc.edu

[16], part III). This makes sense for example for Hammersley's process and exclusion type processes. If the particle system lives in one-dimensional space, we can also look at the system as the height function of an interface, so that the occupation variables of the particle system define the increments of the height function. In terms of the height functions, we start by coupling two systems $\sigma$ and $\zeta$ so that initially $\zeta = \sigma$ to the left of $X(0)$, and $\zeta = \sigma + 1$ to the right of $X(0)$. Then at all times the point $X(t)$ is the boundary of the set $\{x : \zeta(x,t) = \sigma(x,t) + 1\}$.

This last idea generalizes naturally to the multidimensional interface model. We couple two height processes $\sigma$ and $\zeta$ that satisfy $\sigma \leq \zeta \leq \sigma + 1$ at all times, and prove that the boundary of the random set $\{x : \zeta(x,t) = \sigma(x,t) + 1\}$ follows the characteristics of the macroscopic equation.

Laws of large numbers for height functions of asymmetric interface models of the general type considered here have been earlier studied in a handful of papers. A hydrodynamic limit for ballistic deposition was proved in [24], and for models that generalize the exclusion process in [19, 20]. Articles [19, 24] deal with totally asymmetric models, and utilize the path-level variational formulation that generalizes from one-dimensional systems [23]. Article [20] introduces a different approach for partially asymmetric systems. These results are existence results only. In other words convergence to a limiting evolution is shown but nothing explicit about the limit is known, except that it is defined by a Hamilton-Jacobi equation. For partially asymmetric systems in more than one dimension it is presently not even known if the limit is deterministic.

Our motivation for introducing a new model is to have a system for which better results could be proved. An advantage over earlier results is that here we can write down explicitly the macroscopic velocity function up to a constant factor. This is because the process is constructed through a homogeneous space-time Poisson process, so we can simultaneously scale space and time. This is not possible for a lattice model. With an (almost) explicit velocity we can calculate macroscopic profiles, for example see what profiles with shocks and rarefaction fans look like. In the earlier cases at best we know that the velocity function is convex (or concave), but whether the velocity is $C^1$ or strictly convex is a hard open question. Here this question is resolved immediately.

Hammersley's process has been a fruitful model for studying large scale behavior of one-dimensional asymmetric systems, by virtue of its connection with the increasing sequence model. For example, by a combination of the path-level variational construction and the Baik-Deift-Johansson estimates [3], one can presently prove the sharpest out-of-equilibrium fluctuation results for this process [25, 27]. The model introduced in the present paper has a similar connection with a simple combinatorial model, so it may not be too unrealistic to expect some benefit from this in the future.

Before the arrival of the powerful combinatorial and analytic approach pioneered in [3], Hammersley's process was used as a tool for investigating the increasing sequences model. This approach was successful in finding the value of the limiting constant [1, 21] and in large deviations [22]. The proofs relied on explicit knowledge of invariant distributions of Hammersley's process. A similar motivation is possible for us too, and this time the object of interest would be the height of the random partial order. But currently we have no explicit knowledge of steady states of the process introduced here, so we cannot use the process to identify the limiting constant for the random order model.

Recently Cator and Groeneboom [8] developed an approach to the one-dimensional Hammersley process that captures the correct order $t^{1/3}$ of the current fluctua-tions across a characteristic. The argument utilizes precise equilibrium calculations and a time reversal that connects maximizing increasing paths with trajectories of second class particles. In [4] this method is adapted to the totally asymmetric exclusion process. Whether the idea can be applied in multidimensional settings remains to be seen.

*Organization of the paper.* We begin by reminding the reader of the random partial order model, and then proceed to define the process and state the limit theorems. Proofs follow. A technical appendix at the end of the paper shows that the process has a natural state space that is a complete, separable metric space.

## 2. The height of a random partial order

Fix an integer $\nu \geq 2$. Coordinatewise partial orders on $\mathbf{R}^\nu$ are defined for points $x = (x_1, \ldots, x_\nu)$ and $y = (y_1, \ldots, y_\nu)$ by

(1) $\quad x \leq y$ iff $x_i \leq y_i$ for $1 \leq j \leq \nu$, and $x < y$ iff $x_i < y_i$ for $1 \leq j \leq \nu$.

We use interval notation to denote rectangles: $(a, b] = \{x \in \mathbf{R}^\nu : a < x \leq b\}$ for $a < b$ in $\mathbf{R}^\nu$, and similarly for $[a, b]$ and the other types of intervals.

Consider a homogeneous rate 1 Poisson point process in $\mathbf{R}^\nu$. A sequence of Poisson points $p^k$, $1 \leq k \leq m$, is *increasing* if $p^1 < p^2 < \cdots < p^m$ in the coordinatewise sense. For $a < b$ in $\mathbf{R}^\nu$, let $\mathbf{H}(a, b)$ denote the maximal number of Poisson points on an increasing sequence contained in the set $(a, b]$. Let $\mathbf{1} = (1, 1, \ldots, 1) \in \mathbf{R}^\nu$. (This is the only vector we will denote by a boldface.) Kingman's subadditive ergodic theorem and simple moment bounds imply the existence of constants $c_\nu$ such that

(2) $$\lim_{n \to \infty} \frac{1}{n} \mathbf{H}(0, n\mathbf{1}) = c_\nu \qquad \text{a.s.}$$

Presently the only known value is $c_2 = 2$, first proved by Vershik and Kerov [29] and Logan and Shepp [17]. The case $\nu = 2$ is the same as the problem of the longest increasing subsequence of a random permutation, see [2] for a review. Bollobás and Winkler [7] proved that $c_\nu \to e$ as $\nu \to \infty$.

The general case is called the $\nu$-dimensional random partial order, and $\mathbf{H}(0, n\mathbf{1})$ is the *height* of the random partial order. The study of random partial orders was initiated by Winkler [30]. Here is an alternative construction of the random partial order on a fixed (rather than Poisson) number of elements. From the $k!$ linear orders on a set of $k$ elements, choose $\nu$ orders $\prec_1, \prec_2, \ldots, \prec_\nu$ uniformly at random with replacement. Define the random order $\prec$ as the intersection, namely $x \prec y$ iff $x \prec_j y$ for $j = 1, \ldots, \nu$. The height of the random order is the maximal size of a linearly ordered subset. Conditioned on the number $k$ of Poisson points in the cube $(0, n\mathbf{1}]$, $\mathbf{H}(0, n\mathbf{1})$ has the same distribution as the height of the random order $\prec$.

Let us also point out that by the spatial scaling of the Poisson point process, for any $b = (b_1, \ldots, b_\nu) > 0$ in $\mathbf{R}^d$,

(3) $$\lim_{n \to \infty} \frac{1}{n} \mathbf{H}(0, nb) = c_\nu (b_1 b_2 b_3 \cdots b_\nu)^{1/\nu} \qquad \text{a.s.}$$

## 3. The interface process

Fix a spatial dimension $d \geq 2$. Appropriately interpreted, everything we say is true in $d = 1$ also, but does not offer anything significantly new. We describe the evolution of a random, integer-valued height function $\sigma = (\sigma(x))_{x \in \mathbf{R}^d}$. Height values $\pm\infty$ are permitted, so the range of the height function is $\mathbf{Z}^* = \mathbf{Z} \cup \{\pm\infty\}$. The state space of the process is the space $\Sigma$ of functions $\sigma : \mathbf{R}^d \to \mathbf{Z}^*$ that satisfy conditions (i)–(iii):

(4)          (i) Monotonicity: $x \leq y$ in $\mathbf{R}^d$ implies $\sigma(x) \leq \sigma(y)$.

The partial order $x \leq y$ on $\mathbf{R}^d$ is the coordinatewise one defined in Section 2.

(ii) Discontinuities restricted to a locally finite, countable collection of coordinate hyperplanes: for each bounded cube $[-q\mathbf{1}, q\mathbf{1}] \subseteq \mathbf{R}^d$, there are finite partitions

$$-q = s_i^0 < s_i^1 < \cdots < s_i^{m_i} = q$$

along each coordinate direction ($1 \leq i \leq d$), such that any discontinuity point of $\sigma$ in $[-q\mathbf{1}, q\mathbf{1}]$ lies on one of the hyperplanes $\{x \in [-q\mathbf{1}, q\mathbf{1}] : x_i = s_i^k\}$, $1 \leq i \leq d$ and $0 \leq k \leq m_i$.

At discontinuities $\sigma$ is continuous from above: $\sigma(y) \to \sigma(x)$ as $y \to x$ so that $y \geq x$ in $\mathbf{R}^d$. Since $\sigma$ is $\mathbf{Z}^*$-valued, this is the same as saying that $\sigma$ is constant on the left closed, right open rectangles

$$(5) \quad [s^k, s^{k+\mathbf{1}}) \equiv \prod_{i=1}^d [s_i^{k_i}, s_i^{k_i+1}), \quad k = (k_1, k_2, \ldots, k_d) \in \prod_{i=1}^d \{0, 1, 2, \ldots, m_i - 1\},$$

determined by the partitions $\{s_i^k : 0 \leq k \leq m_i\}$, $1 \leq i \leq d$.

(iii) A decay condition "at $-\infty$":

$$(6) \quad \text{for every } b \in \mathbf{R}^d, \quad \lim_{M \to \infty} \sup \left\{ |y|_\infty^{-d/(d+1)} \sigma(y) : y \leq b, |y|_\infty \geq M \right\} = -\infty.$$

The role of the (arbitrary) point $b$ in condition (6) is to confine $y$ so that as the limit is taken, all coordinates of $y$ remain bounded above and at least one of them diverges to $-\infty$. Hence we can think of this as "$y \to -\infty$" in $\mathbf{R}^d$. The $\ell^\infty$ norm on $\mathbf{R}^d$ is $|y|_\infty = \max_{1 \leq i \leq d} |y_i|$.

We can give $\Sigma$ a complete, separable metric. Start with a natural Skorohod metric suggested by condition (ii). On bounded rectangles, this has been considered earlier by Bickel and Wichura [5], among others. This metric is then augmented with sufficient control of the left tail so that convergence in this metric preserves (6). The Borel $\sigma$-field under this metric is generated by the coordinate projections $\sigma \mapsto \sigma(x)$. These matters are discussed in a technical appendix at the end of the paper.

Assume given an initial height function $\sigma \in \Sigma$. To construct the dynamics, assume also given a space-time Poisson point process on $\mathbf{R}^d \times (0, \infty)$. We define the process $\sigma(t) = \{\sigma(x, t) : x \in \mathbf{R}^d\}$ for times $t \in [0, \infty)$ by

$$(7) \quad \sigma(x, t) = \sup_{y : y \leq x} \{\sigma(y) + \mathbf{H}((y, 0), (x, t))\}.$$

The random variable $\mathbf{H}((y, 0), (x, t))$ is the maximal number of Poisson points on an increasing sequence in the space-time rectangle

$$((y, 0), (x, t)] = \{(\eta, s) \in \mathbf{R}^d \times (0, t] : y_i < \eta_i \leq x_i \, (1 \leq i \leq d)\},$$

as defined in Section 2. One can prove that, for almost every realization of the Poisson point process, the supremum in (7) is achieved at some $y$, and $\sigma(t) \in \Sigma$ for all $t > 0$. In particular, if initially $\sigma(x)$ is finite then $\sigma(x, t)$ remains finite for all $0 \leq t < \infty$. And if $\sigma(x) = \pm\infty$, then $\sigma(x,t) = \sigma(x)$ for all $0 \leq t < \infty$. This defines a Markov process on the path space $D([0,\infty),\Sigma)$.

The local effect of the dynamical rule (7) is the following. Suppose $(y,t) \in \mathbf{R}^d \times (0,\infty)$ is a Poisson point, and the state at time $t-$ is $\sigma$. Then at time $t$ the state changes to $\sigma^y$ defined by

(8) $$\sigma^y(x) = \begin{cases} \sigma(x) + 1, & \text{if } x \geq y \text{ and } \sigma(x) = \sigma(y), \\ \sigma(x), & \text{for all other } x \in \mathbf{R}^d. \end{cases}$$

We can express the dynamics succinctly like this: Independently at all $x \in \mathbf{R}^d$, $\sigma(x)$ jumps to $\sigma(x) + 1$ at rate $dx$ ($d$-dimensional volume element). When a jump at $x$ happens, the height function $\sigma$ is updated to $\sigma + 1$ on the set $\{w \in \mathbf{R}^d : w \geq x, \sigma(w) = \sigma(x)\}$ to preserve the monotonicity property (4). It also follows that if $\sigma(y) = \pm\infty$ then $\sigma^y = \sigma$.

We express this in generator language as follows. Suppose $\phi$ is a bounded measurable function on $\Sigma$, and supported on a compact cube $K \subseteq \mathbf{R}^d$. By this we mean that $\phi$ is a measurable function of the coordinates $(\sigma(x))_{x \in K}$. Define the generator $\mathcal{L}$ by

(9) $$\mathcal{L}\phi(\sigma) = \int_{\mathbf{R}^d} [\phi(\sigma^y) - \phi(\sigma)] dy.$$

The next theorem verifies that $\mathcal{L}$ gives the infinitesimal description of the process in one basic sense.

**Theorem 3.1.** *For bounded measurable functions $\phi$ on $\Sigma$, $\sigma \in \Sigma$, and $t > 0$,*

(10) $$E^\sigma[\phi(\sigma(t))] - \phi(\sigma) = \int_0^t E^\sigma[\mathcal{L}\phi(\sigma(s))] ds.$$

*$E^\sigma$ denotes expectation under the path measure $P^\sigma$ of the process defined by (7) and started from state $\sigma$.*

## 4. Hydrodynamic limit for the height process

Let $u_0 : \mathbf{R}^d \to \mathbf{R}$ be a nondecreasing locally Lipschitz continuous function, such that for any $b \in \mathbf{R}^d$,

(11) $$\lim_{M \to \infty} \sup \left\{ |y|_\infty^{-d/(d+1)} u_0(y) : y \leq b, |y|_\infty \geq M \right\} = -\infty.$$

The function $u_0$ represents the initial macroscopic height function. Assume that on some probability space we have a sequence of random initial height functions $\{\sigma_n(y,0) : y \in \mathbf{R}^d\}$, indexed by $n$. Each $\sigma_n(\cdot,0)$ is a.s. an element of the state space $\Sigma$. The sequence satisfies a law of large numbers:

(12) for every $y \in \mathbf{R}^d$, $n^{-1}\sigma_n(ny,0) \to u_0(y)$ as $n \to \infty$, a.s.

Additionally there is the following uniform bound on the decay at $-\infty$:

(13) for every fixed $b \in \mathbf{R}^d$ and $C > 0$, with probability 1 there exist finite, possibly random, $M, N > 0$ such that, if $n \geq N$, $y \leq b$ and $|y|_\infty \geq M$, then $\sigma_n(ny,0) \leq -Cn|y|_\infty^{d/(d+1)}$.

Augment the probability space of the initial $\sigma_n(\cdot, 0)$ by a space-time Poisson point process, and define the processes $\sigma_n(x, t)$ by (7). For $x = (x_1, \ldots, x_d) \geq 0$ in $\mathbf{R}^d$, define
$$g(x) = c_{d+1}(x_1 x_2 x_3 \cdots x_d)^{1/(d+1)}.$$

The constant $c_{d+1}$ is the one from (2), and it comes from the partial order among Poisson points in $d+1$ dimensional space-time rectangles.

Define a function $u(x, t)$ on $\mathbf{R}^d \times [0, \infty)$ by $u(x, 0) = u_0(x)$ and for $t > 0$,

(14) $$u(x, t) = \sup_{y : y \leq x} \{u_0(y) + tg((x - y)/t)\}.$$

The function $u$ is nondecreasing in $x$, increasing in $t$, and locally Lipschitz in $\mathbf{R}^d \times (0, \infty)$.

**Theorem 4.1.** *Suppose $u_0$ is a locally Lipschitz function on $\mathbf{R}^d$ that satisfies (11). Define $u(x, t)$ through (14). Assume that the initial random interfaces $\{\sigma_n(y, 0)\}$ satisfy (12) and (13). Then for all $(x, t) \in \mathbf{R}^d \times [0, \infty)$,*

(15) $$\lim_{n \to \infty} n^{-1} \sigma_n(nx, nt) = u(x, t) \qquad a.s.$$

By the monotonicity of the random height and the continuity of the limiting function, the limit (15) holds simultaneously for all $(x, t)$ outside a single exceptional null event.

Extend $g$ to an u.s.c. concave function on all of $\mathbf{R}^d$ by setting $g \equiv -\infty$ outside $[0, \infty)^d$. Define the constant

(16) $$\kappa_d = \left(\frac{c_{d+1}}{d+1}\right)^{d+1}.$$

The concave conjugate of $g$ is $g^*(\rho) = \inf_x \{x \cdot \rho - g(x)\}$, $\rho \in \mathbf{R}^d$. Let $f = -g^*$. Then $f(\rho) = \infty$ for $\rho \notin (0, \infty)^d$, and

(17) $$f(\rho) = \kappa_d(\rho_1 \rho_2 \cdots \rho_d)^{-1} \qquad \text{for } \rho > 0 \text{ in } \mathbf{R}^d.$$

The Hopf-Lax formula (14) implies that $u$ solves the Hamilton-Jacobi equation (see [10])

(18) $$\partial_t u - f(\nabla u) = 0, \qquad u|_{t=0} = u_0.$$

In other words, $f(\nabla u)$ is the upward velocity of the interface, determined by the local slope.

The most basic case of the hydrodynamic limit starts with $\sigma(y, 0) = 0$ for $y \geq 0$ and $\sigma(y, 0) = -\infty$ otherwise. Then $\sigma(x, t) = \mathbf{H}((0, 0), (x, t))$ for $x \geq 0$ and $-\infty$ otherwise. The limit is $u(x, t) = tg(x/t)$.

## 5. The defect boundary limit

Our objective is to generalize the notion of a second class particle from the one-dimensional context. The particle interpretation does not make sense now. But a second class particle also represents a defect in an interface, and is sometimes called a 'defect tracer.' This point of view we adopt. Given an initial height function $\sigma(y, 0)$, perturb it by increasing the height to $\sigma(y, 0) + 1$ for points $y$ in some set $A(0)$. The boundary of the set $A(0)$ corresponds to a second class particle, so we call it the defect boundary. How does the perturbation set $A(\cdot)$ evolve in time? To describe the behavior of this set under hydrodynamic scaling, we need to look at how the Hamilton-Jacobi equation (18) carries information in time.

For $(x, t) \in \mathbf{R}^d \times (0, \infty)$, let $I(x, t)$ be the set of maximizers in (14):

$$(19) \qquad I(x, t) = \{y \in \mathbf{R}^d : y \leq x,\ u(x, t) = u_0(y) + tg((x - y)/t)\}.$$

Continuity and hypothesis (11) guarantee that $I(x, t)$ is a nonempty compact set. It turns out that these three statements (i)–(iii) are equivalent for a point $(x, t)$: (i) the gradient $\nabla u$ in the $x$-variable exists at $(x, t)$, (ii) $u$ is differentiable at $(x, t)$, and (iii) $I(x, t)$ is a singleton. We call a point $(x, t)$ with $t > 0$ a *shock* if $I(x, t)$ has more than one point.

For $y \in \mathbf{R}^d$ let $W(y, t)$ be the set of points $x \in \mathbf{R}^d$ for which $y$ is a maximizer in the Hopf-Lax formula (14) at time $t$:

$$(20) \qquad W(y, t) = \{x \in \mathbf{R}^d : x \geq y,\ u(x, t) = u_0(y) + tg((x - y)/t)\},$$

and for any subset $B \subseteq \mathbf{R}^d$,

$$(21) \qquad W(B, t) = \bigcup_{y \in B} W(y, t).$$

Given a closed set $B \subseteq \mathbf{R}^d$, let

$$(22) \qquad X(B, t) = W(B, t) \cap W(\overline{B^c}, t).$$

$W(B, t)$ and $W(\overline{B^c}, t)$ are both closed sets. We can characterize $x \in X(B, t)$ as follows: if $(x, t)$ is not a shock then the unique maximizer $\{y\} = I(x, t)$ in (14) lies on the boundary of $B$, while if $(x, t)$ is a shock then $I(x, t)$ intersects both $B$ and $\overline{B^c}$.

If dimension $d = 1$ and $B = [a, \infty) \subseteq \mathbf{R}$, an infinite interval, then $X(B, t)$ is precisely the set of points $x$ for which there exists a forward characteristic $x(\cdot)$ such that $x(0) = a$ and $x(t) = x$. By a forward characteristic we mean a Filippov solution of $dx/dt = f'(\nabla u(x, t))$ [9, 18]. A corresponding characterization of $X(B, t)$ in multiple dimensions does not seem to exist at the moment.

The open $\varepsilon$-neighborhood of a set $B \subseteq \mathbf{R}^d$ is denoted by

$$(23) \qquad B^{(\varepsilon)} = \{x : d(x, y) < \varepsilon \text{ for some } y \in B\}.$$

The distance $d(x, y)$ can be the standard Euclidean distance or another equivalent metric, it makes no difference. Let us write $B^{(-\varepsilon)}$ for the set of $x \in B$ that are at least distance $\varepsilon > 0$ away from the boundary:

$$(24) \qquad B^{(-\varepsilon)} = \{x \in B : d(x, y) \geq \varepsilon \text{ for all } y \notin B\} = \left[(B^c)^{(\varepsilon)}\right]^c.$$

The topological boundary of a closed set $B$ is $\mathrm{bd}B = B \cap \overline{B^c}$.

Suppose two height processes $\sigma(t)$ and $\zeta(t)$ are coupled through the space-time Poisson point process. This means that on some probability space are defined the initial height functions $\sigma(y, 0)$ and $\zeta(y, 0)$, and a space-time Poisson point process which defines all the random variables $\mathbf{H}((y, 0), (x, t))$. Process $\sigma(x, t)$ is defined by (7), and process $\zeta(x, t)$ by the same formula with $\sigma$ replaced by $\zeta$, but with the *same realization* of the variables $\mathbf{H}((y,0),(x,t))$. If initially $\sigma \leq \zeta \leq \sigma + h$ for some constant $h$, then the evolution preserves these inequalities. We can follow the evolution of the "defect set" $A(t)$, defined as $A(t) = \{x : \zeta(x,t) = \sigma(x,t) + h\}$ for $t \geq 0$. This type of a setting we now study in the hydrodynamic context. In the introduction we only discussed the case $h = 1$, but the proof works for general finite $h$.

Now precise assumptions. On some probability space are defined two sequences of initial height functions $\sigma_n(y,0)$ and $\zeta_n(y,0)$. The $\{\sigma_n(y,0)\}$ satisfy the hypotheses (12) and (13) of Theorem 4.1. For some fixed positive integer $h$,

(25) $\qquad \sigma_n(y,0) \leq \zeta_n(y,0) \leq \sigma_n(y,0) + h \quad \text{for all } n \text{ and } y \in \mathbf{R}^d.$

Construct the processes $\sigma_n(t)$ and $\zeta_n(t)$ with the same realizations of the space-time Poisson point process. Then

(26) $\qquad \sigma_n(x,t) \leq \zeta_n(x,t) \leq \sigma_n(x,t) + h \quad \text{for all } n \text{ and } (x,t).$

In particular, $\zeta_n$ and $\sigma_n$ satisfy the same hydrodynamic limit.

Let

(27) $\qquad A_n(t) = \{x \in \mathbf{R}^d : \zeta_n(x,t) = \sigma_n(x,t) + h\}.$

Our objective is to follow the evolution of the set $A_n(t)$ and its boundary $\mathrm{bd}\{A_n(t)\}$. We need an initial assumption at time $t = 0$. Fix a deterministic closed set $B \subseteq \mathbf{R}^d$. We assume that for large $n$, $n^{-1}A_n(0)$ approximates $B$ locally, in the following sense: almost surely, for every compact $K \subseteq \mathbf{R}^d$ and $\varepsilon > 0$,

(28) $\qquad B^{(-\varepsilon)} \cap K \subseteq \{n^{-1}A_n(0)\} \cap K \subseteq B^{(\varepsilon)} \cap K \qquad \text{for all large enough } n.$

**Theorem 5.1.** *Let again $u_0$ satisfy (11) and the processes $\sigma_n$ satisfy (12) and (13) at time zero. Fix a positive integer $h$ and a closed set $B \subseteq \mathbf{R}^d$. Assume that the processes $\sigma_n$ are coupled with processes $\zeta_n$ through a common space-time Poisson point process so that (26) holds. Define $A_n(t)$ by (27) and assume $A_n(0)$ satisfies (28).*

*If $W(B,t) = \emptyset$, then almost surely, for every compact $K \subseteq \mathbf{R}^d$, $A_n(nt) \cap nK = \emptyset$ for all large enough $n$.*

*Suppose $W(B,t) \neq \emptyset$. Then almost surely, for every compact $K \subseteq \mathbf{R}^d$ and $\varepsilon > 0$,*

(29) $\qquad \mathrm{bd}\,\{n^{-1}A_n(nt)\} \cap K \subseteq X(B,t)^{(\varepsilon)} \cap K \qquad \text{for all large enough } n.$

*In addition, suppose no point of $W(\overline{B^c},t)$ is an interior point of $W(B,t)$. Then almost surely, for every compact $K \subseteq \mathbf{R}^d$ and $\varepsilon > 0$,*

(30) $\qquad \begin{aligned} W(B,t)^{(-\varepsilon)} \cap K &\subseteq \{n^{-1}A_n(nt)\} \cap K \\ &\subseteq W(B,t)^{(\varepsilon)} \cap K \qquad \text{for all large enough } n. \end{aligned}$

The additional hypothesis for (30), that no point of $W(\overline{B^c},t)$ is an interior point of $W(B,t)$, prevents $B$ and $\overline{B^c}$ from becoming too entangled at later times. For example, it prohibits the existence of a point $y \in \mathrm{bd}\,B$ such that $W(y,t)$ has nonempty interior ("a rarefaction fan with interior").

## 6. Examples and technical comments

### 6.1. Second class particle analogy

Consider a one-dimensional Hammersley process $z(t) = (z_i(t))_{i \in \mathbf{Z}}$ with labeled particle locations $\cdots \leq z_{-1}(t) \leq z_0(t) \leq z_1(t) \leq \cdots$ on $\mathbf{R}$. In terms of labeled particles, the infinitesimal jump rule is this: $z_i$ jumps to the left at exponential rate $z_i - z_{i-1}$, and when it jumps, its new position $z_i'$ is chosen uniformly at random from the interval $(z_{i-1}, z_i)$. The height function is defined by $\sigma(x, t) = \sup\{i : z_i(t) \leq x\}$ for $x \in \mathbf{R}$.

Now consider another Hammersley process $\tilde{z}(t)$ constructed with the same realization of the space-time Poisson point process as $z(t)$. Assume that at time 0, $\tilde{z}(0)$ has exactly the same particle locations as $z(0)$, plus $h$ additional particles. Then at all later times $\tilde{z}(t)$ will have $h$ particles more than $z(t)$, and relative to the $z(t)$-process, these extra particles behave like second class particles.

Suppose the labeling of the particles is such that $\tilde{z}_i(t) = z_i(t)$ to the left of all the second class particles. Let $X_1(t) \leq \cdots \leq X_h(t)$ be the locations of the second class particles. Then the height functions satisfy $\tilde{\sigma}(x, t) = \sigma(x, t)$ for $x < X_1(t)$, and $\tilde{\sigma}(x, t) = \sigma(x, t) + h$ for $x \geq X_h(t)$. So in this one-dimensional second class particle picture, the set $A(t)$ is the interval $[X_h(t), \infty)$. It has been proved, in the context of one-dimensional asymmetric exclusion, $K$-exclusion and zero-range processes, that in the hydrodynamic limit a second-class particle converges to a characteristic or shock of the macroscopic p.d.e. [12, 18, 26].

Despite this analogy, good properties of the one-dimensional situation are readily lost as we move to higher dimensions. For example, we can begin with a set $A(0)$ that is monotone in the sense that $x \in A(0)$ implies $y \in A(0)$ for all $y \geq x$. But this property can be immediately lost: Suppose a jump happens at $w$ such that $\zeta(w, 0) = \sigma(w, 0)$ but the set $V = \{x \geq w : \sigma(x, 0) = \sigma(w, 0)\}$ intersects $A(0) = \{x : \zeta(x, 0) = \sigma(x, 0) + 1\}$. Then after this event $\zeta = \sigma$ on $V$, and cutting $V$ away from $A(0)$ may have broken its monotonicity.

### 6.2. Examples of the limit in Theorem 5.1

We consider here the simplest macroscopic profiles for which we can explicitly calculate the evolution $W(B, t)$ of a set $B$, and thereby we know the limit of $n^{-1} A_n(nt)$ in Theorem 5.1. These are the flat profile with constant slope, and the cases of shocks and rarefaction fans that have two different slopes. Recall the slope-dependent velocity $f(\rho) = \kappa_d (\rho_1 \rho_2 \rho_3 \cdots \rho_d)^{-1}$ for $\rho \in (0, \infty)^d$, where $\kappa_d$ is the (unknown) constant defined by (2) and (16).

For the second class particle in one-dimensional asymmetric exclusion, these cases were studied in [12, 13].

*Flat profile.* Fix a vector $\rho \in (0, \infty)^d$, and consider the initial profile $u_0(x) = \rho \cdot x$. Then $u(x, t) = \rho \cdot x + t f(\rho)$, for each $(x, t)$ there is a unique maximizer $y(x, t) = x + t \nabla f(\rho)$ in the Hopf-Lax formula, and consequently for any set $B$, $W(B, t) = -t \nabla f(\rho) + B$.

*Shock profile.* Fix two vectors $\lambda, \rho \in (0, \infty)^d$, and let

$$u_0(x) = \begin{cases} \rho \cdot x, & (\rho - \lambda) \cdot x \geq 0, \\ \lambda \cdot x, & (\rho - \lambda) \cdot x \leq 0. \end{cases} \tag{31}$$

Then at later times we have

$$u(x,t) = \begin{cases} \rho \cdot x + tf(\rho), & (\rho - \lambda) \cdot x \geq t(f(\lambda) - f(\rho)), \\ \lambda \cdot x + tf(\lambda), & (\rho - \lambda) \cdot x \leq t(f(\lambda) - f(\rho)). \end{cases}$$

The Hopf-Lax formula is maximized by

$$y = \begin{cases} x + t\nabla f(\rho), & \text{if } (\rho - \lambda) \cdot x \geq t(f(\lambda) - f(\rho)), \\ x + t\nabla f(\lambda), & \text{if } (\rho - \lambda) \cdot x \leq t(f(\lambda) - f(\rho)). \end{cases}$$

In particular, points $(x,t)$ on the hyperplane $(\rho - \lambda) \cdot x = t(f(\lambda) - f(\rho))$ are shocks, and for them both alternatives above are maximizers. In the forward evolution, $W(y,t)$ is either a singleton or empty:

$$W(y,t) = \begin{cases} y - t\nabla f(\rho), \\ \qquad \text{if } (\rho - \lambda) \cdot y \geq t(f(\lambda) - f(\rho)) + t(\rho - \lambda) \cdot \nabla f(\rho), \\ \emptyset, \text{ if } t(f(\lambda) - f(\rho)) + t(\rho - \lambda) \cdot \nabla f(\lambda) < (\rho - \lambda) \cdot y \\ \qquad\qquad\qquad < t(f(\lambda) - f(\rho)) + t(\rho - \lambda) \cdot \nabla f(\rho), \\ y - t\nabla f(\lambda), \\ \qquad \text{if } (\rho - \lambda) \cdot y \leq t(f(\lambda) - f(\rho)) + t(\rho - \lambda) \cdot \nabla f(\lambda). \end{cases}$$

In this situation Theorem 5.1 is valid for all sets $B$.

*Rarefaction fan profile.* Fix two vectors $\lambda, \rho \in (0, \infty)^d$, and let

$$u_0(x) = \begin{cases} \lambda \cdot x, & (\rho - \lambda) \cdot x \geq 0, \\ \rho \cdot x, & (\rho - \lambda) \cdot x \leq 0. \end{cases}$$

For $(x,t)$ such that

$$-t(\rho - \lambda) \cdot \nabla f(\rho) < (\rho - \lambda) \cdot x < -t(\rho - \lambda) \cdot \nabla f(\lambda)$$

there exists a unique $s = s(x,t) \in (0,1)$ such that

$$(\rho - \lambda) \cdot x = -t(\rho - \lambda) \cdot \nabla f(s\lambda + (1-s)\rho).$$

Then at later times the profile can be expressed as

$$u(x,t) = \begin{cases} \rho \cdot x + tf(\rho), & \text{if } (\rho - \lambda) \cdot x \leq -t(\rho - \lambda) \cdot \nabla f(\rho), \\ (s\lambda + (1-s)\rho) \cdot x + tf(s\lambda + (1-s)\rho), & \text{if} \\ \qquad -t(\rho - \lambda) \cdot \nabla f(\rho) < (\rho - \lambda) \cdot x < -t(\rho - \lambda) \cdot \nabla f(\lambda), \\ \lambda \cdot x + tf(\lambda), & \text{if } (\rho - \lambda) \cdot x \geq -t(\rho - \lambda) \cdot \nabla f(\lambda). \end{cases}$$

The forward evolution manifests the rarefaction fan: points $y$ on the hyperplane $(\rho - \lambda) \cdot y = 0$ have $W(y,t)$ given by a curve, while for other points $y$ $W(y,t)$ is a singleton:

$$W(y,t) = \begin{cases} y - t\nabla f(\rho), & \text{if } (\rho - \lambda) \cdot y < 0, \\ \{y - t\nabla f(s\lambda + (1-s)\rho) : 0 \leq s \leq 1\}, & \text{if } (\rho - \lambda) \cdot y = 0, \\ y - t\nabla f(\lambda), & \text{if } (\rho - \lambda) \cdot y > 0. \end{cases}$$

In Theorem 5.1, consider the half-space $B = \{x : (\rho - \lambda) \cdot x \geq 0\}$. Then

$$X(B,t) = \{x : -t(\rho - \lambda) \cdot \nabla f(\rho) \leq (\rho - \lambda) \cdot x \leq -t(\rho - \lambda) \cdot \nabla f(\lambda)\},$$

the "rarefaction strip" in space. Statement (30) is not valid for $B$, because the interior of $X(B,t)$ lies in the interiors of both $W(B,t)$ and $W(\overline{B^c},t)$. Statement (29) is valid, and says that the boundary of $n^{-1}A_n(nt)$ is locally contained in any neighborhood of $X(B,t)$.

In the corresponding one-dimensional setting, Ferrari and Kipnis [13] proved that on the macroscopic scale, the second class particle is uniformly distributed in the rarefaction fan. Their proof depended on explicit calculations with Bernoulli distributions, so presently we cannot approach such precise knowledge of $\mathrm{bd}\{n^{-1}A_n(nt)\}$.

### 6.3. Some random initial conditions

We give here some natural examples of random initial conditions for Theorems 4.1 and 5.1 for the case $d = 2$. We construct these examples from space-time evolutions of one-dimensional Hammersley's process. The space-time coordinates $(y,t)$ of the 1-dimensional process will equal the 2-dimensional spatial coordinates $x = (x_1, x_2)$ of a height function.

*Flat profiles.* In one dimension, Aldous and Diaconis [1] denoted the Hammersley process by $N(y,t)$. The function $y \mapsto N(y,t)$ ($y \in \mathbf{R}$) can be regarded as the counting function of a point process on $\mathbf{R}$. Homogeneous Poisson point processes are invariant for this process.

To construct all flat initial profiles $u_0(x) = \rho \cdot x$ on $\mathbf{R}^2$, we need two parameters that can be adjusted. The rate $\mu$ of the spatial equilibrium of $N(y,t)$ gives one parameter. Another parameter $\tau$ is the jump rate, in other words the rate of the space-time Poisson point process in the graphical construction of $N(y,t)$. Let now $N(y,t)$ be a process in equilibrium, defined for $-\infty < t < \infty$, normalized so that $N(0,0) = 0$, with jump rate $\tau$, and so that the spatial distribution at each fixed time is a homogeneous Poisson process at rate $\mu$. Then the process of particles jumping past a fixed point in space is Poisson at rate $\tau/\mu$ [1], Lemma 8. Consequently $EN(y,t) = \mu y + (\tau/\mu)t$.

This way we can construct a random initial profile whose mean is a given flat initial profile: given $\rho = (\rho_1, \rho_2) \in (0,\infty)^2$, take an equilibrium process $\{N(y,t) : y \in \mathbf{R}, t \in \mathbf{R}\}$ with $\mu = \rho_1$ and $\tau = \rho_1\rho_2$, and define the initial height function for $x = (x_1, x_2) \in \mathbf{R}^2$ by $\sigma((x_1, x_2), 0) = N(x_1, x_2)$.

*Shock profiles.* Next we construct a class of initial shock profiles. Suppose $\rho = (\rho_1, \rho_2)$ and $\lambda = (\lambda_1, \lambda_2)$ satisfy $\rho > \lambda$ and $\rho_1/\rho_2 < \lambda_1/\lambda_2$. Start by constructing the equilibrium Hammersley system $\{N(y,t) : y \in \mathbf{R}, t \in \mathbf{R}\}$ with spatial density $\mu = \lambda_1$ and jump rate $\tau = \lambda_1\lambda_2$. Set $a = (\rho_1 - \lambda_1)/(\rho_2 - \lambda_2) > 0$. Stop each Hammersley particle the first time it hits the space-time line $t = -ay$, and "erase" the entire evolution of $N(y,t)$ above this line. The assumption $\rho_1/\rho_2 < \lambda_1/\lambda_2$ guarantees that each particle eventually hits this line. Now we have constructed the slope-$\lambda$ height function $\sigma((x_1, x_2), 0) = N(x_1, x_2)$ below the line $(\rho - \lambda) \cdot x = 0 \iff x_2 = -ax_1$. (Slope-$\lambda$ in the sense that $E\sigma(x,0) = \lambda \cdot x$.)

To continue the construction, put a rate $\tau' = \rho_1\rho_2$ space-time Poisson point process above the line $t = -ay$ in the space-time picture of the 1-dim Hammersley process. Let the Hammersley particles evolve from their stopped locations on the line $t = -ay$, according to the usual graphical construction [1] of the process, using the rate $\tau'$ space-time Poisson points. The construction is well defined, because given any finite time $T$, $N(y,T)$ is already constructed for $y \leq -T/a$, and for $y > -T/a$ the particle trajectories can be constructed one at a time from left to right, starting with the leftmost particle stopped at a point $(y, -ay)$ for $y > -T/a$.

One can check that defining $\sigma((x_1, x_2), 0) = N(x_1, x_2)$ for $x_2 > -ax_1$ gives the slope-$\rho$ height function above the line $(\rho - \lambda) \cdot x = 0$. Now we have a random initial height function $\sigma(x, 0)$ with mean $E\sigma(x, 0) = u_0(x)$ as in (31).

Finally, we describe a way to define initial configurations for the coupled processes $\zeta$ and $\sigma$ in the context of this shock example. We shall do it so that the set $\{x : \zeta(x, 0) = \sigma(x, 0) + 1\}$ lies inside $B = \{x : x_2 \geq -ax_1\}$, and approximates it closely. Let $\zeta(x, 0)$ be the height function defined above in terms of the $N(y, t)$ constructed in two steps, first below and then above the line $t = -ay$. Let $z_k(t)$ be the trajectories of the labeled Hammersley particles. These trajectories are the level curves of $\zeta(x, 0)$, namely $\zeta((x_1, x_2), 0) \geq k$ iff $z_k(x_2) \leq x_1$. The construction performed above has the property that each $z_k(t)$ crosses the line $t = -ay$ exactly once (the particles were stopped upon first hitting this line, and then continued entirely above the line).

Define new trajectories $z'_k(t)$ as follows: $z'_k(t) = z_k(t)$ below the line $t = -ay$. From the line $t = -ay$ the trajectory $z'_k(t)$ proceeds vertically upward (in the $t$-direction) until it hits the trajectory of $z_{k+1}(t)$. From that point onwards $z'_k(t)$ follows the trajectory of $z_{k+1}(t)$. This is done for all $k$. Let $N'(y, t)$ be the counting function defined by $N'(y, t) = \sup\{k : z'_k(t) \leq y\}$. And then set $\sigma((x_1, x_2), 0) = N'(x_1, x_2)$

The initial height functions $\sigma(x, 0)$ and $\zeta(x, 0)$ thus defined have these properties: $\sigma(x, 0) = \zeta(x, 0)$ for $x_2 \leq -ax_1$. For any point $(x_1, x_2)$ such that $x_2 > -ax_1$ and some particle trajectory $z_k(t)$ passes between $(x_1, -ax_1)$ and $(x_1, x_2)$, $\zeta(x, 0) = \sigma(x, 0) + 1$. This construction satisfies the hypotheses of Theorem 5.1.

### 6.4. Some properties of the multidimensional Hamilton-Jacobi equation

Let $u(x, t)$ be the viscosity solution of the equation $u_t = f(\nabla u)$, defined by the Hopf-Lax formula (14). By assumption, the initial profile $u_0$ is locally Lipschitz and satisfies the decay estimate (11). Hypothesis (11) is tailored to this particular velocity function, and needs to be changed if $f$ is changed.

Part (b) of this lemma will be needed in the proof of Thm. 5.1.

**Lemma 6.1.** (a) *For any compact* $K \subseteq \mathbf{R}^d$, $\bigcup_{x \in K} I(x, t)$ *is compact.*
(b) $W(B, t)$ *is closed for any closed set* $B \subseteq \mathbf{R}^d$.

*Proof.* (a) By (11), as $y \to -\infty$ for $y \leq x$, $u_0(y) + tg((x - y)/t)$ tends to $-\infty$ uniformly over $x$ in a bounded set. Also, the condition inside (19) is preserved by limits because all the functions are continuous. (b) If $W(B, t) \ni x_j \to x$, then by (a) any sequence of maximizers $y_j \in I(x_j, t) \cap B$ has a convergent subsequence. $\square$

The association of $I(x, t)$ to $x$ is not as well-behaved as in one dimension. For example, not only is there no monotonicity, but a simple example can have $x_1 < x_2$ with maximizers $y_i \in I(x_i, t)$ such that $y_2 < y_1$. The local Lipschitz condition on $u_0$ guarantees that each $y \in I(x, t)$ satisfies $y < x$ (i.e. strict inequality for all coordinates).

Properties that are not hard to check include the following. Part (a) of the lemma implies that $u(x, t)$ is locally Lipschitz on $\mathbf{R}^d \times (0, \infty)$. Lipschitz continuity does not necessarily hold down to $t = 0$, but continuity does. $u$ is differentiable at $(x, t)$ iff $I(x, t)$ is a singleton $\{y\}$, and then $\nabla u(x, t) = \nabla g((x - y)/t)$. Also, $\nabla u$ is continuous on the set where it is defined because whenever $(x_n, t_n) \to (x, t)$ and $y_n \in I(x_n, t_n)$, the sequence $\{y_n\}$ is bounded and all limit points lie in $I(x, t)$.

A converse question is when $W(y,t)$ has more than one point. As in one dimension, one can give a criterion based on the regularity of $u_0$ at $y$. The subdifferential $D^-u_0(x)$ and superdifferential $D^+u_0(x)$ of $u_0$ at $x$ are defined by

$$D^-u_0(x) = \left\{ q \in \mathbf{R}^d : \liminf_{y \to x} \frac{u_0(y) - u_0(x) - q \cdot (y-x)}{\|y-x\|} \geq 0 \right\}$$

and

$$D^+u_0(x) = \left\{ p \in \mathbf{R}^d : \limsup_{y \to x} \frac{u_0(y) - u_0(x) - p \cdot (y-x)}{\|y-x\|} \leq 0 \right\}.$$

It is a fact that both $D^\pm u_0(x)$ are nonempty iff $u_0$ is differentiable at $x$, and then $D^\pm u_0(x) = \{\nabla u_0(x)\}$.

One can check that $W(y,t) \subseteq y - t\nabla f(D^+u_0(y))$. Consequently if $D^-u_0(y)$ is nonempty, $W(y,t)$ cannot have more than 1 point. Another fact from one-dimensional systems that also holds in multiple dimensions is that if we restart the evolution at time $s > 0$, then all forward sets $W(y,t)$ are empty or singletons. In other words, if $\tilde{u}$ is a solution with initial profile $\tilde{u}_0$, and we define $u_0(x) = \tilde{u}(x,s)$ and $u(x,t) = \tilde{u}(x,s+t)$, then $D^-u_0(y)$ is never empty. This is because $\nabla g((x-y)/s)$ lies in $D^-\{\tilde{u}(\cdot,s)\}(x)$ for every $y$ that maximizes the Hopf-Lax formula for $\tilde{u}(x,s)$.

## 7. Proof of the generator relation

In this section we prove Theorem 3.1. Throughout the proofs we use the abbreviation

$$x! = x_1 x_2 x_3 \cdots x_d$$

for a point $x = (x_1, \ldots, x_d) \in \mathbf{R}^d$. We make the following definition related to the dynamics of the process. For a height function $\sigma \in \Sigma$ and a point $x \in \mathbf{R}^d$, let

$$(32) \qquad S_x(\sigma) = \begin{cases} \{y \in \mathbf{R}^d : y \leq x, \sigma(y) = \sigma(x)\}, & \text{if } \sigma(x) \text{ is finite,} \\ \emptyset, & \text{if } \sigma(x) = \pm\infty. \end{cases}$$

$S_x(\sigma)$ is the set in space where a Poisson point must arrive in the next instant in order to increase the height value at $x$. Consequently the Lebesgue measure (volume) $|S_x(\sigma)|$ is the instantaneous rate at which the height $\sigma(x)$ jumps up by 1. Since values $\sigma(x) = \pm\infty$ are not changed by the dynamics, it is sensible to set $S_x(\sigma)$ empty in this case. For a set $K$ in $\mathbf{R}^d$ we define

$$(33) \qquad S_K(\sigma) = \bigcup_{x \in K} S_x(\sigma),$$

the set in space where an instantaneous Poisson arrival would change the function $\sigma$ in the set $K$.

We begin with a simple estimate.

**Lemma 7.1.** *Let $x > 0$ in $\mathbf{R}^d$, $t > 0$, and $k$ a positive integer. Then*

$$P\{\mathbf{H}((0,0),(x,t)) \geq k\} \leq \frac{(x!t)^k}{(k!)^{(d+1)}} \leq e^{-k(d+1)}$$

*where the second inequality is valid if $k \geq e^2(x!t)^{1/(d+1)}$. Note that above $(0,0)$ means the space-time point $(0,0) \in \mathbf{R}^d \times [0,\infty)$.*

*Proof.* Let $\gamma = x!t = x_1 x_2 x_3 \cdots x_d \cdot t > 0$ be the volume of the space-time rectangle $(0, x] \times (0, t]$. $k$ uniform points in this $(d+1)$-dimensional rectangle form an increasing chain with probability $(k!)^{-d}$. Thus

$$P\{\mathbf{H}((0,0),(x,t)) \geq k\} \leq \sum_{j: j \geq k} \frac{e^{-\gamma}\gamma^j}{j!}\binom{j}{k}(k!)^{-d} = \gamma^k(k!)^{-(d+1)}$$

$$\leq \gamma^k(k/e)^{-k(d+1)} \leq e^{-k(d+1)}$$

if $k \geq e^2 \gamma^{1/(d+1)}$. $\square$

We need to make a number of definitions that enable us to control the height functions $\sigma \in \Sigma$. For $b \in \mathbf{R}^d$ and $h \in \mathbf{Z}$, let $y^{b,h}(\sigma)$ be the maximal point $y \leq b$ such that the rectangle $[y, b]$ contains the set $\{x \leq b : \sigma(x) \geq h\}$, with $y^{b,h}(\sigma) = b$ if this set is empty. Note that if $y^{b,h}(\sigma) \neq b$ then there exists $x \leq b$ such that $\sigma(x) \geq h$ and $|b - x|_\infty = |b - y^{b,h}(\sigma)|_\infty$.

Throughout this section we consider compact cubes of the type

$$K = [-q\mathbf{1}, q\mathbf{1}] \subseteq \mathbf{R}^d$$

for a fixed number $q > 0$. When the context is clear we may abbreviate $y^h = y^{q\mathbf{1},h}(\sigma)$. Define

$$\lambda_k(\sigma) = \sup_{-\infty < h \leq k-2} (q\mathbf{1} - y^h)! \cdot (k - h)^{-(d+1)}.$$

Property (6) of the state space guarantees that $\lambda_k(\sigma) < \infty$.

The minimal and maximal *finite* height values in $K$ are defined by

$$I(K, \sigma) = \min\{\sigma(x) : x \in K, -\infty < \sigma(x) < \infty\}$$

and

$$J(K, \sigma) = \max\{\sigma(x) : x \in K, -\infty < \sigma(x) < \infty\}.$$

If $\sigma = \pm\infty$ on all of $K$ we interpret $I(K, \sigma) = \infty = -J(K, \sigma)$. Otherwise these quantities are finite because $\sigma$ can take only finitely many values in $K$. If $\sigma$ is finite on all of $K$ then by monotonicity $I(K, \sigma) = \sigma(-q\mathbf{1})$ and $J(K, \sigma) = \sigma(q\mathbf{1})$. Set

$$(34) \qquad \psi_K(\sigma) = \left( \sum_{k=I(K,\sigma)+1}^{J(K,\sigma)+1} \lambda_k^2(\sigma) \right)^{1/2}.$$

If $\sigma = \pm\infty$ on all of $K$ then $\psi_K(\sigma) = 0$.

The next two lemmas are preliminary and illustrate how $\psi_K(\sigma)$ appears as a bound.

**Lemma 7.2.** *For a cube $K = [-q\mathbf{1}, q\mathbf{1}]$ and $\sigma \in \Sigma$,*

$$(35) \qquad |S_K(\sigma)| \leq 2^{d+1}\psi_K(\sigma).$$

*Proof.* If $\sigma = \pm\infty$ on all of $K$ then both sides of (35) are zero by the definitions. Suppose $I(K, \sigma)$ is finite (this is the complementary case). If $x \in K$ and $y \in S_x(\sigma)$, then $y \leq x \leq q\mathbf{1}$ and $\sigma(y) = \sigma(x) \geq I(K, \sigma)$, and consequently $y \in [y^{I(K,\sigma)}, q\mathbf{1}]$. This is true for an arbitrary point $y \in S_K(\sigma)$. Since $y^{I(K,\sigma)-1} \leq y^{I(K,\sigma)}$, we can weaken the conclusion to $S_K(\sigma) \subseteq [y^{I(K,\sigma)-1}, q\mathbf{1}]$ to get

$$|S_K(\sigma)| \leq (q\mathbf{1} - y^{I(K,\sigma)-1})! = \frac{2^{d+1}(q\mathbf{1} - y^{I(K,\sigma)-1})!}{(I(K,\sigma) + 1 - (I(K,\sigma) - 1))^{d+1}}$$

$$\leq 2^{d+1}\lambda_{I(K,\sigma)+1}(\sigma) \leq 2^{d+1}\psi_K(\sigma). \qquad \square$$

**Lemma 7.3.** *Define the event*

(36)
$$G = \{\text{there exist } x \in K \text{ and } y \in \mathbf{R}^d \text{ such that } y < x,$$
$$-\infty < \sigma(y) \leq \sigma(x) - 1 < \infty,$$
$$\text{and } \mathbf{H}((y,0),(x,t)) \geq \sigma(x) + 1 - \sigma(y)\}.$$

*Then for $0 < t < 1/(2e^{d+1}\psi_K(\sigma))$,*

$$P^\sigma(G) \leq 2e^{2(d+1)}\psi_K^2(\sigma)t^2.$$

*Proof.* Let the index $k$ run through the finite values of $\sigma(x)+1$ in $K$, and $h$ represent $\sigma(y)$. Then

$$P^\sigma(G) \leq \sum_{k=I(K,\sigma)+1}^{J(K,\sigma)+1} \sum_{h \leq k-2} P\{\mathbf{H}((y^h,0),(q\mathbf{1},t)) \geq k - h\}.$$

By Lemma 7.1 and the inequality $j! \geq (j/e)^j$, for a fixed $k$ the inner sum becomes

$$\sum_{h \leq k-2} \frac{((q\mathbf{1} - y^h)!t)^{k-h}}{((k-h)!)^{d+1}} \leq \sum_{j \geq 2} \frac{(\lambda_k(\sigma)t)^j j^{(d+1)j}}{(j!)^{d+1}} \leq 2(e^{d+1}\lambda_k(\sigma)t)^2.$$

The assumption on $t$ was used to sum the geometric series. Now sum over $k$. □

We get the first bound on the evolution.

**Lemma 7.4.** *Let $\phi$ be a bounded measurable function on $\Sigma$, supported on a compact cube $K = [-q\mathbf{1}, q\mathbf{1}] \subseteq \mathbf{R}^d$. This means that $\phi(\sigma)$ depends on $\sigma$ only through $(\sigma(x))_{x \in K}$. Then there is a finite constant $C = C(\|\phi\|_\infty)$ such that, for all $\sigma \in \Sigma$ and $t > 0$, the quantity*

$$\Delta_t(\sigma) = E^\sigma[\phi(\sigma(t))] - \phi(\sigma) - t\mathcal{L}\phi(\sigma)$$

*satisfies the bound $|\Delta_t(\sigma)| \leq Ct^2\psi_K^2(\sigma)$, provided $0 \leq t < 1/(2e^{d+1}\psi_K(\sigma))$.*

*Proof.* We may assume that $\sigma$ is not $\pm\infty$ on all of $K$. For otherwise from the definitions

$$E^\sigma[\phi(\sigma(t))] = \phi(\sigma) \quad \text{and} \quad \mathcal{L}\phi(\sigma) = 0,$$

and the lemma is trivially satisfied.

Observe that on the complement $G^c$ of the event defined in (36),

$$\sigma(x,t) = \sup_{y \in S_x(\sigma) \cup \{x\}} \{\sigma(y) + \mathbf{H}((y,0),(x,t))\}$$

for all $x \in K$. (The singleton $\{x\}$ is added to $S_x(\sigma)$ only to accommodate those $x$ for which $\sigma(x) = \pm\infty$ and $S_x(\sigma)$ was defined to be empty.) Consequently on the event $G^c$ the value $\phi(\sigma(t))$ is determined by $\sigma$ and the Poisson points in the space-time region $S_K(\sigma) \times (0,t]$. Let $D_j$ be the event that $S_K(\sigma) \times (0,t]$ contains $j$ space-time Poisson points, and $\overline{D}_2 = (D_0 \cup D_1)^c$ the event that this set contains at least 2 Poisson points. On the event $D_1$, let $Y \in \mathbf{R}^d$ denote the space coordinate of the unique Poisson point, uniformly distributed on $S_K(\sigma)$. Then

$$E^\sigma[\phi(\sigma(t))] = \phi(\sigma) \cdot P^\sigma(G^c \cap D_0) + E[\phi(\sigma^Y) \cdot \mathbf{1}_{G^c \cap D_1}] + O\left(P^\sigma(G) + P^\sigma(\overline{D}_2)\right)$$
$$= \phi(\sigma) + t\mathcal{L}\phi(\sigma) + t^2 \cdot O\left(\psi_K^2(\sigma) + |S_K(\sigma)|^2\right).$$

To get the second equality above, use

$$P^\sigma(D_j) = (j!)^{-1}(t|S_K(\sigma)|)^j \exp(-t|S_K(\sigma)|),$$

Lemma 7.3 for bounding $P^\sigma(G)$, and hide the constant from Lemma 7.3 and $\|\phi\|_\infty$ in the $O$-terms. Proof of the lemma is completed by (35). □

We insert here an intermediate bound on the height $\mathbf{H}$. It is a consequence of Lemma 7.1 and a discretization of space.

**Lemma 7.5.** *Fix $t > 0$, $\alpha \in (0, 1/2)$, and $\beta > e^2 t^{1/(d+1)}$. Then there are finite positive constants $\theta_0$, $C_1$ and $C_2$ such that, for any $\theta \geq \theta_0$,*

$$\begin{aligned}(37) \quad & P\big\{\text{there exist } y < x \text{ such that } |y|_\infty \geq \theta, \ |x - y|_\infty \geq \alpha |y|_\infty, \\ & \text{and } \mathbf{H}((y, 0), (x, t)) \geq \beta |x - y|_\infty^{d/(d+1)}\big\} \leq C_1 \exp(-C_2 \theta^{d/(d+1)}).\end{aligned}$$

For positive $m$, define

$$(38) \qquad \sigma^m(x, t) \equiv \sup_{\substack{y \leq x \\ |y|_\infty \leq m}} \{\sigma(y) + H((y, 0), (x, t))\}$$

**Corollary 7.6.** *Fix a compact cube $K \subseteq \mathbf{R}^d$, $0 < T < \infty$, and initial state $\sigma \in \Sigma$. Then there exists a finite random variable $M$ such that, almost surely, $\sigma(x, t) = \sigma^M(x, t)$ for $(x, t) \in K \times [0, T]$.*

*Proof.* If $\sigma(x) = \pm\infty$ then $y = x$ is the only maximizer needed in the variational formula (7). Thus we may assume that $I(K, \sigma)$ is finite.

Fix $\alpha \in (0, 1/2)$ and $\beta > e^2 T^{1/(d+1)}$. By the boundedness of $K$, (37), and Borel-Cantelli there is a finite random $M$ such that

$$\mathbf{H}((y, 0), (x, T)) \leq \beta |x - y|_\infty^{d/(d+1)}$$

whenever $x \in K$ and $|y|_\infty \geq M$. Increase $M$ further so that $M \geq 1 + |x|$ for all $x \in K$, and $\sigma(y) \leq -(2\beta + |I(K, \sigma)| + 1)|y|_\infty^{-d/(d+1)}$ for all $y$ such that $y \leq x$ for some $x \in K$ and $|y|_\infty \geq M$.

Now suppose $y \leq x$, $x \in K$, $\sigma(x)$ is finite, and $|y|_\infty \geq M$. Then

$$\begin{aligned}\sigma(y) + \mathbf{H}((y, 0), (x, t)) &\leq -(2\beta + |I(K, \sigma)| + 1)|y|_\infty^{-d/(d+1)} + \beta|x - y|_\infty^{d/(d+1)} \\ &\leq -|I(K, \sigma)| - 1 \leq \sigma(x) - 1 \leq \sigma(x, t) - 1.\end{aligned}$$

We see that $y$ cannot participate in the supremum in (7) for any $(x, t) \in K \times [0, T]$. □

To derive the generator formula we need to control the error in Lemma 7.4 uniformly over time, in the form $\Delta_\tau(\sigma(s))$ with $0 \leq s \leq t$ and a small $\tau > 0$. For a fixed $k$, $\lambda_k(\sigma(s))$ is nondecreasing in $s$ because each coordinate of $y^h$ decreases over time. For $q > 0$ and $k \in \mathbf{Z}$ introduce the function

$$(39) \qquad \Psi_{q,k}(\sigma) = \sup_{x \leq q\mathbf{1}} \frac{|x|_\infty^d}{\big(1 \vee \{k - \sigma(x)\}\big)^{d+1}}.$$

A calculation that begins with

$$\lambda_k(\sigma) \leq q^d/2 + \sup_{\substack{h \leq k-2 \\ y^{q\mathbf{1}, h} \neq q\mathbf{1}}} \frac{2^d |y^{q\mathbf{1}, h}|_\infty^d}{(k - h)^{d+1}}$$

shows that
$$\lambda_k(\sigma) \leq q^d/2 + 2^d \Psi_{q,k}(\sigma).$$

Interface heights $\sigma(x,s)$ never decrease with time, and $\Psi_{q,k}(\sigma)$ is nonincreasing in $k$ but nondecreasing in $\sigma$. Therefore we can bound as follows, uniformly over $s \in [0,t]$:

$$\psi_K^2(\sigma(s)) = \sum_{k=I(K,\sigma(s))+1}^{J(K,\sigma(s))+1} \lambda_k^2(\sigma(s))$$
$$\leq \big(J(K,\sigma(t)) - I(K,\sigma(0)) + 1\big) \cdot \max_{I(K,\sigma(0))+1 \leq k \leq J(K,\sigma(t))+1} \lambda_k^2(\sigma(s))$$
(40) $$\leq \big(J(K,\sigma(t)) - I(K,\sigma(0)) + 1\big)^2 (q^{2d} + 1) + 2^{4d} \Psi_{q,I(K,\sigma(0))}^4(\sigma(t)).$$

Above we used the inequality $c(a+b)^2 \leq 2ca^2 + c^2 + b^4$ for $a,b,c \geq 0$. The next lemma implies that the moments $E^\sigma[\Psi_{q,I(K,\sigma)}^p(\sigma(t))]$ are finite for all $p < \infty$.

**Lemma 7.7.** *Let $\sigma$ be an element of the state space $\Sigma$. Fix $t > 0$ and a point $q\mathbf{1} \in \mathbf{R}_+^d$. Then there exists a finite number $v_0(\sigma)$ such that, for $v \geq v_0(\sigma)$,*

(41) $$P^\sigma\big\{\Psi_{q,I(K,\sigma)}(\sigma(t)) > v\big\} \leq C_1 \exp(-C_2 v^{1/(d+1)}),$$

*where the finite positive constants $C_1, C_2$ are the same as in Lemma 7.5 above.*

*Proof.* Choose $\alpha, \beta$ so that (37) is valid. Let
$$\beta_1 = 2\beta + \beta(2\alpha)^{d/(d+1)} + 2.$$

Fix $v_0 = v_0(\sigma) > 0$ so that these requirements are met: $v_0 \geq 1 + |I(K,\sigma)|^{d+1}$, and for all $y \leq x \leq q\mathbf{1}$ such that $|x|_\infty^d \geq v_0$,
$$\sigma(x) \leq -\beta_1 |x|_\infty^{d/(d+1)} \quad \text{and} \quad |y|_\infty \geq |x|_\infty \geq \theta_0.$$

Here $\theta_0$ is the constant that appeared in Lemma 7.5, and we used property (6) of the state space $\Sigma$.

Let $v \geq v_0$. We shall show that the event on the left-hand side of (41) is contained in the event in (37) with $\theta = v^{1/d}$. Suppose the event in (41) happens, so that some $x \leq q\mathbf{1}$ satisfies

(42) $$v^{-1/(d+1)} |x|_\infty^{d(d+1)} > I(K,\sigma) - \sigma(x,t) \quad \text{and} \quad |x|_\infty^d \geq v.$$

Note that the above inequality forces $\sigma(x,t) > -\infty$, while the earlier requirement on $v_0$ forces $\sigma(x) < \infty$, and thereby also $\sigma(x,t) < \infty$. Find a maximizer $y \leq x$ so that
$$\sigma(x,t) = \sigma(y) + \mathbf{H}((y,0),(x,t)).$$

Regarding the location of $y$, we have two cases two consider.

*Case 1.* $y \in [x - 2\alpha|x|_\infty \mathbf{1}, x]$. Let $y' = x - 2\alpha|x|_\infty \mathbf{1}$. Then $|x - y'|_\infty \geq \alpha|y'|_\infty$ by virtue of $\alpha \in (0, 1/2)$. Also $y' \leq x$ so the choices made above imply $|y'|_\infty \geq |x|_\infty \geq v^{1/d}$.

$$\mathbf{H}((y',0),(x,t)) \geq \mathbf{H}((y,0),(x,t)) = \sigma(x,t) - \sigma(y)$$
$$> I(K,\sigma) - v^{-1/(d+1)}|x|_\infty^{d/(d+1)} - \sigma(x)$$
$$\geq (\beta_1 - 2)|x|_\infty^{d/(d+1)} \geq \beta(2\alpha)^{d/(d+1)} |x|_\infty^{d/(d+1)}$$
$$= \beta|x - y'|_\infty^{d/(d+1)}.$$

In addition to (42), we used $-v^{-1/(d+1)} \geq -1$, $-\sigma(y) \geq -\sigma(x) \geq \beta_1 |x|_\infty^{d/(d+1)}$, and $I(K,\sigma) \geq -v_0^{1/(d+1)} \geq -|x|_\infty^{d/(d+1)}$.

*Case 2.* $y \notin [x - 2\alpha|x|_\infty \mathbf{1}, x]$. This implies $|x - y|_\infty \geq \alpha |y|_\infty$.

$$\mathbf{H}((y,0),(x,t)) = \sigma(x,t) - \sigma(y) > I(K,\sigma) - v^{-1/(d+1)}|x|_\infty^{d/(d+1)} - \sigma(y)$$
$$\geq -|x|_\infty^{d/(d+1)} - v^{-1/(d+1)}|x|_\infty^{d/(d+1)} + \beta_1|y|_\infty^{d/(d+1)}$$
$$\geq 2\beta|y|_\infty^{d/(d+1)} \geq \beta|x - y|_\infty^{d/(d+1)}.$$

We conclude that the event in (41) lies inside the event in (37) with $\theta = v^{1/d}$, as long as $v \geq v_0$, and the inequality in (41) follows from (37). □

**Corollary 7.8.** *Let $K$ be a compact cube, $\varepsilon > 0$, and $0 < t < \infty$. Then there exists a deterministic compact cube $L$ such that*

$$P^\sigma\{S_K(\sigma(s)) \subseteq L \text{ for all } s \in [0,t]\} \geq 1 - \varepsilon.$$

*Proof.* For $0 \leq s \leq t$, $x \in S_K(\sigma(s))$ implies that $I(K, \sigma(s))$ is finite, $x \leq q\mathbf{1}$ and $\sigma(x,s) \geq I(K, \sigma(s))$. Consequently

$$|x|_\infty^d \leq \Psi_{q,I(K,\sigma(s))}(\sigma(s)) \leq \Psi_{q,I(K,\sigma)}(\sigma(t)).$$

Thus given $\varepsilon$, we can choose $L = [-m\mathbf{1}, m\mathbf{1}]$ with $m$ picked by Lemma 7.7 so that $P^\sigma\{\Psi_{q,I(K,\sigma)}(\sigma(t)) > m^d\} < \varepsilon$. □

We are ready for the last stage of the proof of Theorem 3.1.

**Proposition 7.9.** *Let $\phi$ be a bounded measurable function on $\Sigma$ supported on the compact cube $K = [-q\mathbf{1}, q\mathbf{1}]$ of $\mathbf{R}^d$, and $\sigma \in \Sigma$. Then*

$$(43) \qquad E^\sigma[\phi(\sigma(t))] - \phi(\sigma) = \int_0^t E^\sigma[\mathcal{L}\phi(\sigma(s))]ds.$$

*Proof.* Pick a small $\tau > 0$ so that $t = m\tau$ for an integer $m$, and denote the partition by $s_j = j\tau$. By the Markov property,

$$E^\sigma[\phi(\sigma(t))] - \phi(\sigma) = E^\sigma\left[\sum_{j=0}^{m-1}\left\{E^{\sigma(s_j)}[\phi(\sigma(\tau))] - \phi(\sigma(s_j))\right\}\right]$$
$$= E^\sigma\left[\int_0^t \sum_{j=0}^{m-1} \mathbf{1}_{(s_j, s_{j+1}]}(s)\mathcal{L}\phi(\sigma(s_{j+1}))ds\right]$$
$$(44) \qquad + \tau(\mathcal{L}\phi(\sigma) - E^\sigma[\mathcal{L}\phi(\sigma(t))]) + E^\sigma\left[\sum_{j=0}^{m-1}\Delta_\tau(\sigma(s_j))\right],$$

where the terms $\Delta_\tau(\sigma(s_j))$ are as defined in Lemma 7.4.

We wish to argue that, as $m \to \infty$ and simultaneously $\tau \to 0$, expression (44) after the last equality sign converges to the right-hand side of (43).

Note first that $\mathcal{L}\phi(\sigma)$ is determined by the restriction of $\sigma$ to the set $S_K(\sigma) \cup K$. By Corollary 7.8 there exists a fixed compact set $L$ such that $S_K(\sigma(s)) \cup K \subseteq L$ for $0 \leq s \leq t$ with probability at least $1 - \varepsilon$. By Corollary 7.6, the time evolution $\{\sigma(x,s) : x \in L, 0 \leq s \leq t\}$ is determined by the finitely many Poisson points in the random compact rectangle $[-M, M]^d \times [0, t]$. Consequently the process $\mathcal{L}\phi(\sigma(s))$ is piecewise constant in time, and then the integrand $\sum_{j=0}^{m-1} \mathbf{1}_{(s_j, s_{j+1}]}(s)\mathcal{L}\phi(\sigma(s_{j+1}))$ converges to $\mathcal{L}\phi(\sigma(s))$ pointwise as $m \to \infty$. This happens on an event with probability at least $1 - \varepsilon$, hence almost surely after letting $\varepsilon \to 0$.

To extend the convergence to the expectation and to handle the error terms, we show that

$$(45) \qquad E^\sigma \Big[ \sup_{0 \leq s \leq t} \psi_K^2(\sigma(s)) \Big] < \infty.$$

Before proving (45), let us see why it is sufficient. Since

$$(46) \qquad |\mathcal{L}\phi(\sigma)| \leq 2\|\phi\|_\infty |S_K(\sigma)|,$$

(35) and (45) imply that also the first expectation after the equality sign in (44) converges, by dominated convergence. The second and third terms of (44) vanish, through a combination of Lemma 7.4, (35), and (45).

By the bound in (40) for $\sup_{0 \leq s \leq t} \psi_K^2(\sigma(s))$ and by Lemma 7.7, it only remains to show that

$$E^\sigma \big[ \big( J(K, \sigma(t)) - I(K, \sigma(0)) + 1 \big)^2 \big] < \infty.$$

This follows from property (6) of $\sigma$ and the bounds for $\mathbf{H}$ in Lemmas 7.1 and 7.5. We omit the proof since it is not different in spirit than the estimates we already developed. □

This completes the proof of Theorem 3.1.

## 8. Proof of the limit for the height function

Introduce the scaling into the variational formula (7) and write it as

$$(47) \qquad \sigma_n(nx, nt) = \sup_{y \in \mathbf{R}^d : y \leq x} \{\sigma_n(ny, 0) + \mathbf{H}((ny, 0), (nx, nt))\}.$$

**Lemma 8.1.** *Assume the processes $\sigma_n$ satisfy (12) and (13). Fix a finite $T > 0$ and a point $b \in \mathbf{R}^d$ such that $b > 0$, and consider the bounded rectangle $[-b, b] \subseteq \mathbf{R}^d$. Then with probability 1 there exist a random $N < \infty$ and a random point $a \in \mathbf{R}^d$ such that*

$$(48) \qquad \sigma_n(nx, nt) = \sup_{y \in [a, x]} \{\sigma_n(ny, 0) + \mathbf{H}((ny, 0), (nx, nt))\}$$

*for $x \in [-b, b]$, $t \in (0, T]$, $n \geq N$.*

*Proof.* For $\beta \geq e^2 T^{1/(d+1)}$ and $b \in \mathbf{R}^d$ fixed, one can deduce from Lemma 7.1 and Borel-Cantelli that, almost surely, for large enough $n$,

$$\mathbf{H}((ni, 0), (nb, nt)) \leq \beta n|b - i|_\infty^{d/(d+1)}$$

for all $i \in \mathbf{Z}^d$ such that $i \leq b$ and $|i - b|_\infty \geq 1$. If $y \in \mathbf{R}^d$ satisfies $y \leq b$ and $|y - b|_\infty \geq 1$, we can take $i = [y]$ (coordinatewise integer parts of $y$) and see that

$$(49) \qquad \mathbf{H}((ny, 0), (nb, nt)) \leq \beta n + \beta n|b - y|_\infty^{d/(d+1)}$$

for all such $y$.

In assumption (13) choose $C > \beta$ so that $-C+(2+|b|_\infty^{d/(d+1)})\beta < u_0(-b)-1$. Let $N$ and $M$ be as given by (13), but increase $M$ further to guarantee $M \geq 1$. Now take $a \in \mathbf{R}^d$ far enough below $-b$ so that, if $y \leq b$ but $y \geq a$ fails, then $|y|_\infty \geq M$. [Since assumption (13) permits a random $M > 0$, here we may need to choose a random $a \in \mathbf{R}^d$.] Then by (13), if $y \leq b$ but $y \geq a$ fails, then $\sigma_n(ny,0) \leq -Cn|y|_\infty^{d/(d+1)}$.

Now suppose $x \in [-b,b]$, $y \leq x$, but $y \geq a$ fails. Then

$$\begin{aligned}
\sigma_n(ny,0) + \mathbf{H}((ny,0),(nx,nt)) & \\
&\leq \sigma_n(ny,0) + \mathbf{H}((ny,0),(nb,nt)) \\
&\leq -Cn|y|_\infty^{d/(d+1)} + \beta n + \beta n|b-y|_\infty^{d/(d+1)} \\
&\leq n\left((-C+\beta)|y|_\infty^{d/(d+1)} + \beta + \beta|b|_\infty^{d/(d+1)}\right) \\
&\leq nu_0(-b) - n \leq \sigma_n(-nb,0) - n/2 \\
&\qquad\qquad\text{[by assumption (12), for large enough } n\text{]} \\
&\leq \sigma_n(nx,0) - n/2 \qquad \text{[by monotonicity]}.
\end{aligned}$$

This shows that in the variational formula (47) the point $y = x$ strictly dominates all $y$ outside $[a,x]$. □

Starting with (48) the limit (15) is proved (i) by partitioning $[a,x]$ into small rectangles, (ii) by using monotonicity of the random variables, and the monotonicity and continuity of the limit, and (iii) by appealing to the assumed initial limits (12) and to

(50) $\quad n^{-1}\mathbf{H}((ny,0),(nx,nt)) \to c_{d+1}((x-y)!t)^{1/(d+1)} = tg((x-y)/t)$ a.s.

To derive the limit in (50) from (3) one has to fill in a technical step because in (50) the lower left corner of the rectangle $(ny, nx] \times (0, nt]$ moves as $n$ grows. One can argue around this complication in at least two different ways: (a) The Kesten-Hammersley lemma [28], page 20, from subadditive theory gives a.s. convergence along a subsequence, and then one fills in to get the full sequence. This approach was used in [24]. (b) Alternatively, one can use Borel-Cantelli if summable deviation bounds are available. These can be obtained by combining Theorems 3 and 9 from Bollobás and Brightwell [6].

## 9. Proof of the defect boundary limit

In view of the variational equation (7), let us say $\sigma(x,t)$ has a maximizer $y$ if $y \leq x$ and $\sigma(x,t) = \sigma(y,0) + \mathbf{H}((y,0),(x,t))$.

**Lemma 9.1.** *Suppose two processes $\sigma$ and $\zeta$ are coupled through the space-time Poisson point process.*

*(a) For a positive integer $m$, let $D_m(t) = \{x : \zeta(x,t) \geq \sigma(x,t) + m\}$. Then if $x \in D_m(t)$, $\zeta(x,t)$ cannot have a maximizer $y \in D_m(0)^c$. And if $x \in D_m(t)^c$, $\sigma(x,t)$ cannot have a maximizer $y \in D_m(0)$.*

*(b) In particular, suppose initially $\sigma(y,0) \leq \zeta(y,0) \leq \sigma(y,0) + h$ for all $y \in \mathbf{R}^d$, for a fixed positive integer $h$. Then this property is preserved for all time. If we write*

$$A(t) = \{x : \zeta(x,t) = \sigma(x,t) + h\},$$

*then*

(51) $\qquad A(t) = \{x : \sigma(x,t) \text{ has a maximizer } y \in A(0)\}.$

(c) *If $h = 1$ in part (b), we get additionally that*

(52) $$A(t)^c = \{x : \zeta(x,t) \text{ has a maximizer } y \in A(0)^c \}.$$

*Proof.* (a) Suppose $x \in D_m(t)$, $y \in D_m(0)^c$, and $\zeta(x,t) = \zeta(y,0) + \mathbf{H}((y,0),(x,t))$. Then by the definition of $D_m(t)$,

$$\sigma(x,t) \leq \zeta(x,t) - m = \zeta(y,0) - m + \mathbf{H}((y,0),(x,t)) \leq \sigma(y,0) + \mathbf{H}((y,0),(x,t)) - 1$$

which contradicts the variational equation (7). Thus $\zeta(x,t)$ cannot have a maximizer $y \in D_m(0)^c$. The other part of (a) is proved similarly.

(b) Monotonicity implies that $\sigma(x,t) \leq \zeta(x,t) \leq \sigma(x,t) + h$ for all $(x,t)$, so $A(t) = D_h(t)$. Suppose $x \in A(t)$. By (a) $\zeta(x,t)$ cannot have a maximizer $y \in A(0)^c$, and so $\zeta(x,t)$ has a maximizer $y \in A(0)$. Consequently

$$\sigma(x,t) = \zeta(x,t) - h = \zeta(y,0) - h + \mathbf{H}((y,0),(x,t)) = \sigma(y,0) + \mathbf{H}((y,0),(x,t)),$$

which says that $\sigma(x,t)$ has a maximizer $y \in A(0)$. On the other hand, if $\sigma(x,t)$ has a maximizer $y \in A(0)$, then by (a) again $x \notin A(t)^c$. This proves (51).

(c) Now $A(t) = D_1(t)$ and $A(t)^c = \{x : \sigma(x,t) = \zeta(x,t)\}$. If $\zeta(x,t)$ has a maximizer $y \in A(0)^c$, then by part (a) $x \notin A(t)$. While if $x \in A(t)^c$, again by part (a) $\sigma(x,t)$ must have a maximizer $y \in A(0)^c$, which then also is a maximizer for $\zeta(x,t)$. This proves (52). $\square$

Assume the sequence of processes $\sigma_n(\cdot)$ satisfies the hypotheses of the hydrodynamic limit Theorem 4.1 which we proved in Section 8. The defect set $A_n(t)$ was defined through the $(\sigma_n, \zeta_n)$ coupling by (27). By (51) above, we can equivalently define it by

(53) $$A_n(t) = \{x : \sigma_n(x,t) \text{ has a maximizer } y \in A_n(0) \}.$$

In the next lemma we take the point of view that *some* sequence of sets that depend on $\omega$ has been defined by (53), and ignore the $(\sigma_n, \zeta_n)$ coupling definition.

**Lemma 9.2.** *Let $B \subseteq \mathbf{R}^d$ be a closed set. Suppose that for almost every sample point $\omega$ in the underlying probability space, a sequence of sets $A_n(0) = A_n(0; \omega)$ is defined, and has this property: for every compact $K \subseteq \mathbf{R}^d$ and $\varepsilon > 0$,*

(54) $$\{n^{-1} A_n(0)\} \cap K \subseteq B^{(\varepsilon)} \cap K \quad \text{for all large enough } n.$$

*Suppose the sets $A_n(t)$ satisfy (53) and fix $t > 0$. Then almost surely, for every compact $K \subseteq \mathbf{R}^d$ and $\varepsilon > 0$,*

(55) $$\{n^{-1} A_n(nt)\} \cap K \subseteq W(B,t)^{(\varepsilon)} \cap K \quad \text{for all large enough } n.$$

*In particular, if $W(B,t) = \emptyset$, then (55) implies that $\{n^{-1} A_n(nt)\} \cap K = \emptyset$ for all large enough $n$.*

*Proof.* Fix a sample point $\omega$ such that assumption (54) is valid, the conclusion of Lemma 8.1 is valid for all $b \in \mathbf{Z}_+^d$, and we have the limits

(56) $$n^{-1} \sigma_n(nx, nt) \to u(x,t) \quad \text{for all } (x,t),$$

and

(57) $\quad n^{-1}\mathbf{H}((ny,0),(nx,nt)) \to tg((x-y)/t) \qquad$ for all $y,x,t$.

Almost every $\omega$ satisfies these requirements, by the a.s. limits (50) and (15), by monotonicity, and by the continuity of the limiting functions. It suffices to prove (55) for this fixed $\omega$.

To contradict (55), suppose there is a subsequence $n_j$ and points $x_j \in K$ such that $n_j x_j \in A_{n_j}(n_j t)$ but $x_j \notin W(B,t)^{(\varepsilon)}$. Note that this also contradicts $\{n^{-1}A_n(nt)\} \cap K = \emptyset$ in case $W(B,t) = \emptyset$, so the empty set case is also proved by the contradiction we derive.

Let $n_j y_j \in A_{n_j}(0)$ be a maximizer for $\sigma_{n_j}(n_j x_j, n_j t)$. Since the $x_j$'s are bounded, so are the $y_j$'s by Lemma 8.1, and we can pass to a subsequence (again denoted by $\{j\}$) such that the limits $x_j \to x$ and $y_j \to y$ exist. By the assumptions on $x_j$, $x \notin W(B,t)$. For any $\varepsilon > 0$, $y_j \in B^{(\varepsilon)}$ for large enough $j$, so $y \in B$ by the closedness of $B$.

Fix points $x' < x''$ and $y' < y''$ so that $x' < x < x''$ and $y' < y < y''$ in the partial order of $\mathbf{R}^d$. Then for large enough $j$, $x' < x_j < x''$ and $y' < y_j < y''$. By the choice of $y_j$,

$$\sigma_{n_j}(n_j x_j, n_j t) = \sigma_{n_j}(n_j y_j, 0) + \mathbf{H}((n_j y_j, 0), (n_j x_j, n_j t))$$

from which follows, by the monotonicity of the processes,

$$n_j^{-1}\sigma_{n_j}(n_j x', n_j t) \leq n_j^{-1}\sigma_{n_j}(n_j x_j, n_j t)$$
$$\leq n_j^{-1}\sigma_{n_j}(n_j y'', 0) + n_j^{-1}\mathbf{H}((n_j y', 0),(n_j x'', n_j t)).$$

Now let $n_j \to \infty$ and use the limits (56) and (57) to obtain

$$u(x',t) \leq u_0(y'') + tg((x''-y')/t).$$

We may let $x', x'' \to x$ and $y', y'' \to y$, and then by continuity $u(x,t) \leq u_0(y) + tg((x-y)/t)$. This is incompatible with having $x \notin W(B,t)$ and $y \in B$. This contradiction shows that, for the fixed $\omega$, (55) holds. $\square$

We prove statement (29) of Theorem 5.1. The assumption is that

(58) $\quad B^{(-\varepsilon)} \cap K \subseteq \{n^{-1}A_n(0)\} \cap K \subseteq B^{(\varepsilon)} \cap K \qquad$ for all large enough $n$.

We introduce an auxiliary process $\xi_n(x,t)$. Initially set

(59) $\quad \xi_n(y,0) = \begin{cases} \sigma_n(y,0), & y \notin A_n(0) \\ \sigma_n(y,0)+1, & y \in A_n(0). \end{cases}$

$\xi_n(y,0)$ is a well-defined random element of the state space $\Sigma$ because $A_n(0)$ is defined (27) in terms of $\zeta_n(y,0)$ which lies in $\Sigma$. Couple the process $\xi_n$ with $\sigma_n$ and $\zeta_n$ through the common space-time Poisson points. Then

$$\sigma_n(x,t) \leq \xi_n(x,t) \leq \sigma_n(x,t) + 1.$$

By part (b) of Lemma 9.1, $A_n(t)$ that satisfies (53) also satisfies

$$A_n(t) = \{x : \xi_n(x,t) = \sigma_n(x,t) + 1\}.$$

Then by part (c) of Lemma 9.1,

(60) $$A_n(t)^c = \{x : \xi_n(x,t) \text{ has a maximizer } y \in A_n(0)^c \}.$$

The first inclusion of assumption (58) implies that $n^{-1}A_n(0)^c \cap K \subseteq (\overline{B^c})^{(\varepsilon)} \cap K$ for large $n$. The processes $\xi_n$ inherit all the hydrodynamic properties of the processes $\sigma_n$. Thus by (60) we may apply Lemma 9.2 to the sets $A_n(nt)^c$ and the processes $\xi_n(nt)$ to get

(61) $$n^{-1}A_n(nt)^c \cap K \subseteq W(\overline{B^c},t)^{(\delta)} \cap K$$

for large enough $n$. By (55) and (61),

$$\text{bd}\,\{n^{-1}A_n(nt)\} \cap K \subseteq W(B,t)^{(\delta)} \cap W(\overline{B^c},t)^{(\delta)} \cap K$$

for large $n$. For small enough $\delta > 0$, the set on the right is contained in $[W(B,t) \cap W(\overline{B^c},t)]^{(\varepsilon)} \cap K = X(B,t)^{(\varepsilon)} \cap K$. This proves (29).

To complete the proof of Theorem 5.1, it remains to prove

(62)
$$W(B,t)^{(-\varepsilon)} \cap K \subseteq n^{-1}A_n(nt) \cap K \subseteq W(B,t)^{(\varepsilon)} \cap K \qquad \text{for all large enough } n$$

under the further assumption that no point of $W(\overline{B^c},t)$ is an interior point of $W(B,t)$.

The second inclusion of (62) we already obtained in Lemma 9.2. (61) implies

$$\left[W(\overline{B^c},t)^{(\delta)}\right]^c \cap K \subseteq n^{-1}A_n(nt) \cap K.$$

It remains to check that, given $\varepsilon > 0$,

$$W(B,t)^{(-\varepsilon)} \cap K \subseteq \left[W(\overline{B^c},t)^{(\delta)}\right]^c \cap K$$

for sufficiently small $\delta > 0$. Suppose not, so that for a sequence $\delta_j \searrow 0$ there exist $x_j \in W(B,t)^{(-\varepsilon)} \cap W(\overline{B^c},t)^{(\delta_j)} \cap K$. By Lemma 6.1 the set $W(\overline{B^c},t)$ is closed. Hence passing to a convergent subsequence $x_j \to x$ gives a point $x \in W(\overline{B^c},t)$ which is an interior point of $W(B,t)$, contrary to the hypothesis.

## 10. Technical appendix: the state space of the process

We develop the state space in two steps: first describe the multidimensional Skorohod type metric we need, and then amend the metric to provide control over the left tail of the height function. This Skorohod type space has been used earlier (see [5] and their references).

### 10.1. A Skorohod type space in multiple dimensions

Let $(X,r)$ be a complete, separable metric space, with metric $r(x,y) \leq 1$. Let $D = D(\mathbf{R}^d, X)$ denote the space of functions $\sigma : \mathbf{R}^d \to X$ with this property: for every bounded rectangle $[a,b] \subseteq \mathbf{R}^d$ and $\varepsilon > 0$, there exist finite partitions

$$a_i = s_i^0 < s_i^1 < \cdots < s_i^{m_i} = b_i$$

of each coordinate axis ($1 \leq i \leq d$) such that the variation of $\sigma$ in the partition rectangles is at most $\varepsilon$: for each $k = (k_1, k_2, \ldots, k_d) \in \prod_{i=1}^{d}\{0, 1, 2, \ldots, m_i - 1\}$,

$$\sup\{r(\sigma(x), \sigma(y)) : s_i^{k_i} \leq x_i, y_i < s_i^{k_i+1} \ (1 \leq i \leq d)\} \leq \varepsilon. \tag{63}$$

Note that the partition rectangles are closed on the left. This implies that $\sigma$ is continuous from above: $\sigma(y) \to \sigma(x)$ as $y \to x$ in $\mathbf{R}^d$ so that $y \geq x$; and limits exist from strictly below: $\lim \sigma(y)$ exists as $y \to x$ in $\mathbf{R}^d$ so that $y < x$ (strict inequality for each coordinate).

We shall employ this notation for truncation in $\mathbf{R}^d$: for real $u > 0$ and $x = (x_1, \ldots, x_d) \in \mathbf{R}^d$,

$$[x]_u = \big((x_1 \wedge u) \vee (-u), (x_2 \wedge u) \vee (-u), \ldots, (x_d \wedge u) \vee (-u)\big).$$

Let $\Lambda$ be the collection of bijective, strictly increasing Lipschitz functions $\lambda : \mathbf{R}^d \to \mathbf{R}^d$ that satisfy these requirements: $\lambda$ is of the type $\lambda(x_1, \ldots, x_d) = (\lambda_1(x_1), \ldots, \lambda_d(x_d))$ where each $\lambda_i : \mathbf{R} \to \mathbf{R}$ is bijective, strictly increasing and Lipschitz; and

$$\gamma(\lambda) = \gamma_0(\lambda) + \gamma_1(\lambda) < \infty$$

where the quantities $\gamma_0(\lambda)$ and $\gamma_1(\lambda)$ are defined by

$$\gamma_0(\lambda) = \sum_{i=1}^{d} \sup_{s,t \in \mathbf{R}} \left|\log \frac{\lambda_i(t) - \lambda_i(s)}{t - s}\right|$$

and

$$\gamma_1(\lambda) = \int_0^\infty e^{-u}\Big(1 \wedge \sup_{x \in \mathbf{R}^d} \big|[\lambda(x)]_u - [x]_u\big|_\infty\Big)\,du.$$

For $\rho, \sigma \in D$, $\lambda \in \Lambda$ and $u > 0$, define

$$d(\rho, \sigma, \lambda, u) = \sup_{x \in \mathbf{R}^d} r\big(\rho([x]_u), \sigma([\lambda(x)]_u)\big).$$

And then

$$d_S(\rho, \sigma) = \inf_{\lambda \in \Lambda} \left\{\gamma(\lambda) \vee \int_0^\infty e^{-u} d(\rho, \sigma, \lambda, u)\,du\right\}. \tag{64}$$

The definition was arranged so that $\gamma(\lambda^{-1}) = \gamma(\lambda)$ and $\gamma(\lambda \circ \mu) \leq \gamma(\lambda) + \gamma(\mu)$, so the proof in [11], Section 3.5, can be repeated to show that $d_S$ is a metric.

It is clear that if a sequence of functions $\sigma_n$ from $D$ converges to an arbitrary function $\sigma : \mathbf{R}^d \to X$, and this convergence happens uniformly on compact subsets of $\mathbf{R}^d$, then $\sigma \in D$. Furthermore, we also get convergence in the $d_S$-metric, as the next lemma indicates. This lemma is needed in the proof that $(D, d_S)$ is complete.

**Lemma 10.1.** *Suppose $\sigma_n, \sigma \in D$. Then $d_S(\sigma_n, \sigma) \to 0$ iff there exist $\lambda_n \in \Lambda$ such that $\gamma(\lambda_n) \to 0$ and*

$$r\big(\sigma_n(x), \sigma(\lambda_n(x))\big) \to 0$$

*uniformly over $x$ in compact subsets of $\mathbf{R}^d$.*

*Proof.* We prove $d_S(\sigma_n, \sigma) \to 0$ assuming the second condition, and leave the other direction to the reader. For each rectangle $[-M\mathbf{1}, M\mathbf{1})$, $M = 1, 2, 3, \ldots$, and each $\varepsilon = 1/K$, $K = 1, 2, 3, \ldots$, fix the partitions $\{s_i^k\}$ that appear in the definition (63) of $\sigma \in D$. Pick a real $u > 0$ so that neither $u$ nor $-u$ is among these countably many partition points.

$$d(\sigma_n, \sigma, \lambda_n, u) = \sup_{x \in \mathbf{R}^d} r\big(\sigma_n([x]_u), \sigma([\lambda_n(x)]_u)\big)$$
$$\leq \sup_{x \in \mathbf{R}^d} r\big(\sigma_n([x]_u), \sigma(\lambda_n([x]_u))\big) + \sup_{x \in \mathbf{R}^d} r\big(\sigma(\lambda_n([x]_u)), \sigma([\lambda_n(x)]_u)\big).$$

The first term after the inequality vanishes as $n \to \infty$, by assumption.

Let $\varepsilon = 1/K > 0$, pick a large rectangle $[-M\mathbf{1}, M\mathbf{1})$ that contains $[-u\mathbf{1}, u\mathbf{1}]$ well inside its interior, and for this rectangle and $\varepsilon$ pick the finite partitions that satisfy (63) for $\sigma$, and do not contain $\pm u$. Let $\delta > 0$ be such that none of these finitely many partition points lie in $(\pm u - \delta, \pm u + \delta)$. If $n$ is large enough, then $\sup_{x \in [-M\mathbf{1}, M\mathbf{1}]} |\lambda_n(x) - x| < \delta$, and one can check that $\lambda_n([x]_u)$ and $[\lambda_n(x)]_u$ lie in the same partition rectangle, for each $x \in \mathbf{R}^d$. Thus

$$\sup_{x \in \mathbf{R}^d} r\left(\sigma(\lambda_n([x]_u)), \sigma([\lambda_n(x)]_u)\right) \leq \varepsilon.$$

We have shown that $d(\sigma_n, \sigma, \lambda_n, u) \to 0$ for a.e. $u > 0$. □

With this lemma, one can follow the proof in [11], page 121, to show that $(D, d_S)$ is complete. Separability of $(D, d_S)$ would also be easy to prove. Next, we take this Skorohod type space as starting point, and define the state space $\Sigma$ for the height process.

### 10.2. The state space for the height process

In the setting of the previous subsection, take $S = \mathbf{Z}^* = \mathbf{Z} \cup \{\pm\infty\}$ with the discrete metric $r(x, y) = \mathbf{1}\{x \neq y\}$. Let $\Sigma$ be the space of functions $\sigma \in D(\mathbf{R}^d, \mathbf{Z}^*)$ that are nondecreasing $[\sigma(x) \leq \sigma(y)$ if $x \leq y$ in $\mathbf{R}^d]$ and decay to $-\infty$ sufficiently fast at $-\infty$, namely

(65) for every $b \in \mathbf{R}^d$, $\lim_{M \to \infty} \sup \left\{ |y|_\infty^{-d/(d+1)} \sigma(y) : y \leq b, |y|_\infty \geq M \right\} = -\infty$.

Condition (65) is not preserved by convergence in the $d_S$ metric, so we need to fix the metric.

For $\sigma \in \Sigma$, $h \in \mathbf{Z}$, and $b \in \mathbf{R}^d$, let $y^{b,h}(\sigma)$ be the maximal $y \leq b$ in $\mathbf{R}^d$ such that the rectangle $[y, b]$ contains the set $\{x \leq b : \sigma(x) \geq h\}$. Condition (65) guarantees that such a finite $y^{b,h}(\sigma)$ exists. In fact, (65) is equivalent to

(66) for every $b \in \mathbf{R}^d$, $\lim_{h \to -\infty} |h|^{-(d+1)/d} |y^{b,h}(\sigma)|_\infty = 0$.

For $\rho, \sigma \in \Sigma$ and $b \in \mathbf{R}^d$, define

$$\theta_b(\rho, \sigma) = \sup_{h \leq -1} |h|^{-(d+1)/d} \cdot |y^{b,h}(\rho) - y^{b,h}(\sigma)|_\infty$$

and

$$\Theta(\rho, \sigma) = \int_{\mathbf{R}^d} e^{-|b|_\infty} \big(1 \wedge \theta_b(\rho, \sigma)\big) db.$$

$\Theta(\rho, \sigma)$ satisfies the triangle inequality, is symmetric, and $\Theta(\sigma, \sigma) = 0$, so we can define a metric on $\Sigma$ by

$$d_\Sigma(\rho, \sigma) = \Theta(\rho, \sigma) + d_S(\rho, \sigma).$$

The effect of the $\Theta(\rho, \sigma)$ term in the metric is the following.

*Growth model* 229**Lemma 10.2.** *Suppose $d_S(\sigma_n, \sigma) \to 0$. Then $d_\Sigma(\sigma_n, \sigma) \to 0$ iff for every $b \in \mathbf{R}^d$,*

(67) $$\lim_{h \to -\infty} \sup_n |h|^{-(d+1)/d} |y^{b,h}(\sigma_n)|_\infty = 0,$$

*or equivalently, for every $b \in \mathbf{R}^d$*

(68) $$\lim_{M \to \infty} \sup_n \sup_{\substack{y \leq b \\ |y|_\infty \geq M}} \frac{\sigma_n(y)}{|y|_\infty^{d/(d+1)}} = -\infty.$$

We leave the proof of the above lemma to the reader. Lemmas 10.1 and 10.2 together give a natural characterization of convergence in $\Sigma$.

**Lemma 10.3.** *The Borel $\sigma$-field $\mathcal{B}_\Sigma$ is the same as the $\sigma$-field $\mathcal{F}$ generated by the coordinate projections $\sigma \mapsto \sigma(x)$.*

*Proof.* The sets $\{x : \sigma(x) \geq h\}$ are closed, so the functions $\sigma \mapsto \sigma(x)$ are upper semicontinuous. This implies $\mathcal{F} \subseteq \mathcal{B}_\Sigma$.

For the other direction one shows that for a fixed $\rho \in \Sigma$, the function $\sigma \mapsto d_\Sigma(\rho, \sigma)$ is $\mathcal{F}$-measurable. This implies that the balls $\{\sigma \in \Sigma : d_\Sigma(\rho, \sigma) < r\}$ are $\mathcal{F}$-measurable. Once we argue below that $\Sigma$ is separable, this suffices for $\mathcal{B}_\Sigma \subseteq \mathcal{F}$.

To show the $\mathcal{F}$-measurability of $\sigma \mapsto d_S(\rho, \sigma)$ one can adapt the argument from page 128 of [11]. To show the $\mathcal{F}$-measurability of $\sigma \mapsto \Theta(\rho, \sigma)$, one can start by arguing the joint $\mathcal{B}_{\mathbf{R}^d} \otimes \mathcal{F}$-measurability of the map $(b, \sigma) \mapsto y^{b,h}(\sigma)$ from $\mathbf{R}^d \times \Sigma$ into $\mathbf{R}^d$. We leave the details. □

The remaining work is to check that $(\Sigma, d_\Sigma)$ is a complete separable metric space.

**Proposition 10.4.** *The space $(\Sigma, d_\Sigma)$ is complete.*

We prove this proposition in several stages. Let $\{\sigma_n\}$ be a Cauchy sequence in the $d_\Sigma$ metric. By the completeness of $(D, d_S)$, we already know there exists a $\sigma \in D(\mathbf{R}^d, \mathbf{Z}^*)$ such that $d_S(\sigma_n, \sigma) \to 0$. We need to show that (i) $\sigma \in \Sigma$ and (ii) $\Theta(\sigma_n, \sigma) \to 0$.

Following the completeness proof for Skorohod space in [11], page 121, we may extract a subsequence, denoted again by $\sigma_n$, together with a sequence of Lipschitz functions $\psi_n \in \Lambda$ (actually labeled $\mu_n^{-1}$ in [11]), such that

(69) $$\gamma(\psi_n) < 2^{1-n}$$

and

(70) $$\sigma_n(\psi_n(x)) \to \sigma(x) \text{ uniformly on compact sets.}$$

**Step 1.** $\sigma \in \Sigma$.

Fix $b \in \mathbf{R}^d$, for which we shall show (66). It suffices to consider $b > 0$. Let $b^k = b + k\mathbf{1}$. By passing to a further subsequence we may assume $\Theta(\sigma_n, \sigma_{n+1}) < e^{-n^2}$. Fix $n_0$ so that

(71) $$\exp\bigl(|b^2|_\infty + d(n+1) - n^2\bigr) < 2^{-n} \text{ for all } n \geq n_0.$$

**Lemma 10.5.** *For $n \geq n_0$ there exist points $\beta^n$ in $\mathbf{R}^d$ such that $b^1 < \beta^{n+1} < \beta^n < b^2$, and $\theta_{\beta^n}(\sigma_n, \sigma_{n+1}) < 2^{-n}$.*

*Proof.* Let $\alpha_n = b^1 + e^{-n} \cdot \mathbf{1}$ in $\mathbf{R}^d$.

$$e^{-n^2} \geq \Theta(\sigma_n, \sigma_{n+1})$$
$$\geq \inf_{x \in (\alpha_{n+1}, \alpha_n)} \{1 \wedge \theta_x(\sigma_n, \sigma_{n+1})\} \cdot e^{-|b^2|_\infty} \cdot \mathrm{Leb}_d\{x : \alpha_{n+1} < x < \alpha_n\},$$

where

$$\mathrm{Leb}_d\{x : \alpha_{n+1} < x < \alpha_n\} = (e^{-n} - e^{-n-1})^d \geq e^{-d(n+1)}$$

is the $d$-dimensional Lebesgue measure of the open rectangle $(\alpha_{n+1}, \alpha_n)$. This implies it is possible to choose a point $\beta^n \in (\alpha_{n+1}, \alpha_n)$ so that $\theta_{\beta^n}(\sigma_n, \sigma_{n+1}) < 2^{-n}$. □

$\beta^{n+1} < \beta^n$ implies $y^{\beta^{n+1}, h}(\sigma_{n+1}) \geq y^{\beta^n, h}(\sigma_{n+1}) - (\beta^n - \beta^{n+1})$. For each fixed $h \leq -1$, applying the above Lemma inductively gives for $n \geq n_0$:

$$\begin{aligned} y^{\beta^{n+1}, h}(\sigma_{n+1}) &\geq y^{\beta^n, h}(\sigma_{n+1}) - (\beta^n - \beta^{n+1}) \\ &\geq y^{\beta^n, h}(\sigma_n) - |h|^{(d+1)/d} 2^{-n} \cdot \mathbf{1} - (\beta^n - \beta^{n+1}) \\ &\geq \cdots \geq y^{\beta^{n_0}, h}(\sigma_{n_0}) - |h|^{(d+1)/d} \sum_{k=n_0}^{n} 2^{-k} \cdot \mathbf{1} - (\beta^{n_0} - \beta^{n+1}), \end{aligned}$$

from which then

(72) $$\inf_{n \geq n_0} y^{b^1, h}(\sigma_n) \geq y^{b^2, h}(\sigma_{n_0}) - |h|^{(d+1)/d} 2^{1-n_0} \cdot \mathbf{1} - (b^2 - b^1).$$

Now fix $h \leq -1$ for the moment. By (72) we may fix a rectangle $[y^1, b^1]$ that contains the sets $\{x \leq b^1 : \sigma_n(x) \geq h\}$ for all $n \geq n_0$. Let $Q = [y^1 - \mathbf{1}, b^1 + \mathbf{1}]$ be a larger rectangle such that each point in $[y^1, b^1]$ is at least distance 1 from $Q^c$. By (69) and (70) we may pick $n$ large enough so that $|\psi_n(x) - x| < 1/4$ and $\sigma_n(\psi_n(x)) = \sigma(x)$ for $x \in Q$. [Equality because $\mathbf{Z}^*$ has the discrete metric.]

We can now argue that if $x \leq b$ and $\sigma(x) \geq h$, then necessarily $x \in Q$, $\psi_n(x) \leq b^1$ and $\sigma_n(\psi_n(x)) \geq h$, which implies by (72) that

$$x \geq \psi_n(x) - (1/4)\mathbf{1} \geq y^{b^1, h}(\sigma_n) - (1/4)\mathbf{1} \geq y^{b^2, h}(\sigma_{n_0}) - (5/4 + |h|^{(d+1)/d} 2^{1-n_0}) \cdot \mathbf{1}.$$

This can be repeated for each $h \leq -1$, with $n_0$ fixed. Thus for all $h \leq -1$,

$$|y^{b, h}(\sigma)| \leq |b| \vee \left(|y^{b^2, h}(\sigma_{n_0})| + 5/4 + |h|^{(d+1)/d} 2^{1-n_0}\right),$$

and then, since $\sigma_{n_0} \in \Sigma$,

$$\lim_{h \to -\infty} |h|^{-(d+1)/d} |y^{b, h}(\sigma)| \leq 2^{1-n_0}.$$

Since $n_0$ can be taken arbitrarily large, (66) follows for $\sigma$, and thereby $\sigma \in \Sigma$.

**Step 2.** $\Theta(\sigma_n, \sigma) \to 0$.

As for Step 1, let us assume that we have picked a subsequence $\sigma_n$ that satisfies (69) and (70) and $\Theta(\sigma_{n+1}, \sigma_n) < e^{-n^2}$. Let $\phi_n = \psi_n^{-1}$. If we prove $\Theta(\sigma_n, \sigma) \to 0$ along this subsequence, then the Cauchy assumption and triangle inequality give it for the full sequence.

Fix an arbitrary index $n_1$ and a small $0 < \varepsilon_0 < 1$. Fix also $\beta \in \mathbf{R}^d$. For each $h \leq -1$, fix a rectangle $[y^h, \beta]$ that contains the sets $\{x \leq \beta : \sigma_n(x) \geq h\}$ for each $\sigma_n$ for $n \geq n_1$, and also for $\sigma$, which Step 1 just showed lies in $\Sigma$. This can be done for each fixed $h$ because by (72) there exists $n_0 = n_0(\beta)$ defined by (71) so that the points $y^{\beta,h}(\sigma_n)$ are bounded below for $n \geq n_0$. Then if necessary decrease $y^h$ further so that

$$y^h \leq y^{\beta,h}(\sigma_1) \wedge y^{\beta,h}(\sigma_2) \wedge \cdots \wedge y^{\beta,h}(\sigma_{n_0-1}).$$

Let $Q_{h,k} = [y^h - k\mathbf{1}, \beta + k\mathbf{1}]$ be larger rectangles.

On the rectangles $Q_{h,2}$, $h \leq -1$, construct the finite partition for $\sigma$ which satisfies (63) for $\varepsilon = 1/2$, so that the discrete metric forces $\sigma$ to be constant on the partition rectangles. Consider a point $b = (b_1, b_2, \ldots, b_d) < \beta$ with the property that no coordinate of $b$ equals any one of the (countably many) partition points. This restriction excludes only a Lebesgue null set of points $b$.

Find $\varepsilon_1 = \varepsilon_1(\beta, b, h) > 0$ such that the intervals $(b_i - \varepsilon_1, b_i + \varepsilon_1)$ contain none of the finitely many partition points that pertain to the rectangle $Q_{h,2}$. Pick $n = n(\beta, b, h) > n_1$ such that $\sigma_n(\psi_n(x)) = \sigma(x)$ and $|\psi_n(x) - x| < (\varepsilon_0 \wedge \varepsilon_1)/4$ for $x \in Q_{h,2}$. Since the maps $\psi, \phi$ do not carry any points of $[y^h, \beta]$ out of $Q_{h,1}$, $y^{h,b}(\sigma_n) = y^{h,b}(\sigma \circ \phi_n)$. It follows that

$$|y^{h,b}(\sigma) - y^{h,b}(\sigma_n)| = |y^{h,b}(\sigma) - y^{h,b}(\sigma \circ \phi_n)| < \varepsilon_0.$$

The last inequality above is justified as follows: The only way it could fail is that $\sigma$ (or $\sigma \circ \phi_n$) has a point $x \leq b$ with height $\geq h$, and $\sigma \circ \phi_n$ (respectively, $\sigma$) does not. These cannot happen because the maps $\psi, \phi$ cannot carry a partition point from one side of $b_i$ to the other side, along any coordinate direction $i$.

Now we have for a.e. $b < \beta$ and each $h \leq -1$, with $n = n(\beta, b, h) > n_1$:

$$\begin{aligned}
&|h|^{-(d+1)/d}|y^{h,b}(\sigma) - y^{h,b}(\sigma_{n_1})| \\
&\quad \leq |h|^{-(d+1)/d}|y^{h,b}(\sigma) - y^{h,b}(\sigma_n)| + |h|^{-(d+1)/d}|y^{h,b}(\sigma_n) - y^{h,b}(\sigma_{n_1})| \\
&\quad \leq \varepsilon_0 + \theta_b(\sigma_n, \sigma_{n_1}) \\
&\quad \leq \varepsilon_0 + \sup_{m:m>n_1} \theta_b(\sigma_m, \sigma_{n_1}).
\end{aligned}$$

The last line has no more dependence on $\beta$ or $h$. Since $\beta$ was arbitrary, this holds for a.e. $b \in \mathbf{R}^d$. Take supremum over $h \leq -1$ on the left, to get

$$1 \wedge \theta_b(\sigma, \sigma_{n_1}) \leq \varepsilon_0 + \sup_{m:m>n_1} \{1 \wedge \theta_b(\sigma_m, \sigma_{n_1})\} \qquad \text{for a.e. } b.$$

Integrate to get

$$\begin{aligned}
\Theta(\sigma, \sigma_{n_1}) &= \int_{\mathbf{R}^d} e^{-|b|_\infty} \{1 \wedge \theta_b(\sigma, \sigma_{n_1})\} \, db \\
&\leq C\varepsilon_0 + \int_{\mathbf{R}^d} e^{-|b|_\infty} \sup_{m:m>n_1} \{1 \wedge \theta_b(\sigma_m, \sigma_{n_1})\} \, db \\
&\leq C\varepsilon_0 + \int_{\mathbf{R}^d} e^{-|b|_\infty} \sup_{m:m>n_1} \left\{\sum_{k=n_1}^{m-1} 1 \wedge \theta_b(\sigma_{k+1}, \sigma_k)\right\} db \\
&= C\varepsilon_0 + \sum_{k=n_1}^{\infty} \int_{\mathbf{R}^d} e^{-|b|_\infty} \{1 \wedge \theta_b(\sigma_{k+1}, \sigma_k)\} \, db \\
&= C\varepsilon_0 + \sum_{k=n_1}^{\infty} \Theta(\sigma_{k+1}, \sigma_k) \leq C\varepsilon_0 + \sum_{k=n_1}^{\infty} e^{-k^2},
\end{aligned}$$

where $C = \int_{\mathbf{R}^d} e^{-|b|_\infty} db$. Since $n_1$ was an arbitrary index, we have

$$\limsup_{n_1 \to \infty} \Theta(\sigma_{n_1}, \sigma) \leq C\varepsilon_0.$$

Since $\varepsilon_0$ was arbitrary, Step 2 is completed, and Proposition 10.4 thereby proved.

We outline how to construct a countable dense set in $(\Sigma, d_\Sigma)$. Fix $a < b$ in $\mathbf{Z}^d$. In the rectangle $[a,b] \subseteq \mathbf{R}^d$, consider the (countably many) finite rational partitions of each coordinate axis. For each such partition of $[a,b]$ into rectangles, consider all the nondecreasing assignments of values from $\mathbf{Z}^*$ to the rectangles. Extend the functions $\widehat{\sigma}$ thus defined to all of $\mathbf{R}^d$ in some fashion, but so that they are nondecreasing and $\mathbf{Z}^*$-valued. Repeat this for all rectangles $[a,b]$ with integer corners. This gives a countable set $\widehat{D}$ of elements of $D(\mathbf{R}^d, \mathbf{Z}^*)$. Finally, each such $\widehat{\sigma} \in \widehat{D}$ yields countably many elements $\widetilde{\sigma} \in \Sigma$ by setting

$$\widetilde{\sigma}(x) = \begin{cases} -\infty, & \widehat{\sigma}(x) < h \\ \widehat{\sigma}(x), & \widehat{\sigma}(x) \geq h \end{cases}$$

for all $h \in \mathbf{Z}$. All these $\widetilde{\sigma}$ together form a countable set $\widetilde{\Sigma} \subseteq \Sigma$.

Now given an arbitrary $\sigma \in \Sigma$, it can be approximated by an element $\widetilde{\sigma} \in \widetilde{\Sigma}$ arbitrarily closely (in the sense that $\sigma = \widetilde{\sigma} \circ \phi$ for a map $\phi \in \Lambda$ close to the identity) on any given finite rectangle $[-\beta, \beta]$, and so that $y^{b,h}(\widetilde{\sigma})$ is close to $y^{b,h}(\sigma)$ for all $b$ in this rectangle, for any given range $h_0 \leq h \leq -1$. Since $|h|^{-(d+1)/d}|y^{\beta,h}(\sigma)| < \varepsilon$ for $h \leq h_0$ for an appropriately chosen $h_0$, this suffices to make both $d(\sigma, \widetilde{\sigma}, \phi, u)$ and $\theta_b(\sigma, \widetilde{\sigma})$ small for a range of $u > 0$ and $b \in \mathbf{R}^d$. To get close under the metric $d_\Sigma$ it suffices to approximate in a bounded set of $u$'s and $b$'s, so it can be checked that $\widetilde{\Sigma}$ is dense in $\Sigma$.

**Acknowledgments.** The author thanks Tom Kurtz for valuable suggestions and anonymous referees for careful readings of the manuscript. Hermann Rost has also studied the process described here but has not published his results.